\documentclass[12pt,francais]{article}
\usepackage[latin1]{inputenc}
\usepackage{babel}
\usepackage{amsmath,amsthm,amsfonts,amssymb,amscd,epsf}
\usepackage{graphicx}
\usepackage{color}

\evensidemargin = \oddsidemargin
\advance\textwidth by -2in
\makeatletter
\gdef\thmhead@plain#1#2#3{%
  \thmname{#1}\thmnumber{\@ifnotempty{#1}{ }#2}%
  \thmnote{ {\mdseries#3}}}
\let\thmhead\thmhead@plain
\makeatother
\theoremstyle{plain}
\newtheorem*{theoreme*}{Théorème}
\newtheorem{theoreme}{Théorème}[section]
\newtheorem{proposition}[theoreme]{Proposition}
\newtheorem{lemme}[theoreme]{Lemme}
\newtheorem{lemme-d'elagage}[theoreme]{Lemme d'élagage}

\theoremstyle{definition}

\theoremstyle{remark}
\newtheorem{remarque}{Remarque}
\newtheorem{conjecture}{Conjecture}

\newtheorem*{exemple}{Exemple}

\newtheorem*{resume}{Résumé}

\def\ly{\fontencoding{U}\fontfamily{lasy}\fontseries{m}\fontshape{n}\selectfont}
\def\guil#1{\leavevmode\hbox{{\ly(\kern-0.20em(\kern+0.20em}}\nobreak{}\,#1\,%
  \nobreak\hbox{{\ly\kern+0.20em)\kern-0.20em)}}}
\def\Alinea#1{\hfill\break%
  \hbox to \parindent{\hss{\textup{#1}}\enspace}\ignorespaces}
\def\alinea#1{\par\noindent%
  \hbox to \parindent{\hss{\textup{#1}}\enspace}\ignorespaces}

\def\Chi{\setbox0=\hbox{$\chi$} \mathord{\raise\dp0\hbox{$\chi$}}}

\def\N{{\mathbb N}}

\def\Q{{\mathbb Q}}
\def\R{{\mathbb R}}

\def\F{{\mathcal F}}

\def\Z{{\mathbb Z}}
\def\CP^#1{\mathbf P^{#1}(\mathbf C)}
\def\RP^#1{\mathbf P^{#1}(\mathbf R)}

\def\Int{\mathop{\mathrm{Int}}\nolimits}

\def\adresse#1{\def\nl{\egroup\egroup\hbox\bgroup\itshape\bgroup}
  \par \noindent
  \hbox to \textwidth{\hskip .5\textwidth \vbox{\small
  \hbox \bgroup\itshape\bgroup #1 \egroup\egroup}\hfill}}

\title{Paires de structures de contact sur les vari\'et\'es de dimension trois}

\author{Vincent Colin et Sebasti\~ao Firmo}
\date{}

\begin{document}
\maketitle

\begin{resume}
On introduit une notion de paire positive de structures de contact
sur les vari\'et\'es
de dimension trois qui g\'en\'eralise celle de \cite{ET,Mi1,Mi2}. Une telle
paire \guil{normale} donne naissance \`a un champ de plans continu et
localement int\'egrable $\lambda$.
On montre que si $\lambda$ est uniquement int\'egrable et
si les structures de contact sont tendues, alors le feuilletage
int\'egral de $\lambda$ est sans composante de Reeb
d'\^ame homologue \`a z\'ero. De plus, dans ce cas, la vari\'et\'e ambiante porte un
feuilletage sans composante de Reeb.
On d\'emontre \'egalement un th\'eor\`eme de  stabilit\'e \guil{\`a la
Reeb} pour les paires positives de structures tendues.

\smallskip

\noindent
\textit{Mots cl\'es}~:
structure de contact, paire, feuilletage, tendu, composante de Reeb

\noindent
\textit{Codes AMS}~:
57R17, 57M50, 57R30.
\end{resume}

\section{Introduction}
En dimension trois, Eliashberg-Thurston~\cite{ET} et Mitsumatsu~\cite{Mi1,Mi2}
associent \`a toute paire $(\xi_+ ,\xi_- )$ de structures de contact
 transversales, o\`u $\xi_+$ est positive et $\xi_-$ n\'egative, une paire de champs de
plans $(\lambda_+ ,\lambda_- )$
transversaux, continus et localement
(non uniquement) int\'egrables~: par tout point passe un germe
de surface int\'egrale de $\lambda_\pm$. L'intersection
 $\xi_+ \cap \xi_-$ est alors dirig\'ee par un champ de vecteurs
\guil{conform\'ement} Anosov. R\'eciproquement,
cette situation se rencontre
lorsqu'on \'etudie les directions  stables et instables
$\lambda_\pm$
d'un champ de vecteurs Anosov, les structures $\xi_\pm$
apparaissant comme \guil{plans m\'edians} des $\lambda_\pm$.

Du point de vue de la g\'eom\'etrie de contact, le fait d'imposer
\`a $\xi_+$ et $\xi_-$ d'\^etre des structures transversales
est trop contraignant. Le but du pr\'esent article est d'\'etudier
le cas de paires $(\xi_+ ,\xi_- )$ o\`u l'on rel\^ache cette
condition pour la remplacer par une propri\'et\'e
de co\" \i ncidence  positive~:
les champs de plans  $\xi_+$ et $\xi_-$ ne sont jamais \guil{dos-\`a-dos}.
Autrement dit, ils sont transversaux \`a un
m\^eme champ de droites $\mathcal{D}$.
Dans cette situation, on perd l'existence d'un des deux champs int\'egrable
$\lambda_\pm$, mais on en pr\'eserve un, $\lambda$, coinc\'e entre
$\xi_+$ et $\xi_-$ et lui aussi transversal \`a $\mathcal{D}$.
Comme auparavant, le champ $\lambda$
est localement int\'egrable et en g\'en\'eral  seulement continu~: il n'y a pas unicit\'e
locale des surfaces int\'egrales, ni donc {\it a priori} de
feuilletage int\'egral global.
Lorsque cette unicit\'e est av\'er\'ee, par exemple sous les
conditions du th\'eor\`eme~\ref{theoreme : transitif}, les propri\'et\'es de
rigidit\'e des structures $\xi_+$ et $\xi_-$ rejaillissent sur celle
du feuilletage int\'egral~: si celles-l\`a sont tendues, celui-ci
est sans composante de Reeb d'\^ame homologue \`a z\'ero et la vari\'et\'e ambiante
porte un feuilletage sans composante de Reeb
(th\'eor\`eme~\ref{theoreme : reeb}).
Notre leitmotiv est que cette notion de paire  positive
de structures de contact tendues pourrait
\^etre une bonne g\'en\'eralisation de celle
de feuilletage sans composante de Reeb.
En particulier, d'apr\`es Eliashberg et Thurston~\cite{ET},
tout feuilletage tendu (et m\^eme tout feuilletage sans composante
de Reeb~\cite{Co2}) est limite de paires positives de structures de contact
tendues.
Pour conclure, on d\'emontre pour les paires positives de structures de contact
tendues un analogue du th\'eor\`eme de stabilit\'e de Reeb
pour les feuilletages (th\'eor\`eme~\ref{t:sphere}).
Ce r\'esultat semble faire
\'echo \`a  la th\'eorie des courbes holomorphes. Son pendant
a en effet \'et\'e d\'emontr\'e par Eliashberg et Hofer~\cite{EH},
{\it via} une m\'ethode de remplissage par des disques holomorphes.
Par certains aspects, les r\'esultats de cet article peuvent
\^etre vus comme une version topologique des feuilletages
d'\'energie finie de Hofer, Wyzocki et Zehnder~\cite{HWZ}.
La notion de paire positive est \'egalement bien adapt\'ee
\`a celle  de surface branch\'ee, connue dans le monde
int\'egrable pour rendre compte des propri\'et\'es des
laminations. Zannad tire  b\'en\'efice
de ce parall\`ele dans \cite{Za1,Za2} en trouvant une condition suffisante,
inspir\'ee de la g\'eom\'etrie de contact, pour qu'une
surface branch\'ee porte une lamination.

\medskip

{\it Remerciements.} Ce travail a \'et\'e rendu possible par l'accord
France-Br\'esil, gr\^ace auquel nous avons pu s\'ejourner
\`a l'Universit\'e Federal Fluminense et \`a l'Universit\'e
de Nantes. Il  a \'egalement b\'en\'efici\'e
du support de l'Institut Universitaire de France et de l'ANR Symplexe.
Nous souhaitons saluer ces institutions
pour leur soutien. Le rapporteur anonyme d'une premi\`ere version de ce texte
y a relev\'e de nombreuses erreurs. Cette nouvelle mouture doit beaucoup
\`a la qualit\'e de son travail. Nous lui adressons nos plus vifs
remerciements.

\section{Paires de structures de contact}

\subsection{Structures de contact}
Une {\it structure de contact} (orientable) sur une vari\'et\'e orient\'ee $V$ de
dimension trois est un champ de plans orientable lisse $\xi$, noyau d'une $1$-forme
$\alpha$ dont le produit ext\'erieur avec $d\alpha$ ne s'annule pas.
Le signe de $\xi$ est celui  de $\alpha \wedge d\alpha$,
rapport\'e \`a l'orientation de $V$.
\`A l'oppos\'e des feuilletages, les structures de contact supportent
bien les d\'eformations~: la condition de contact est ouverte
pour la topologie $C^1$ et, d'apr\`es un th\'eor\`eme
de Gray~\cite{Gr}, tout chemin de structures de contact
est le fait d'une isotopie de $V$ issue de l'identit\'e.

Une structure de contact $\xi$ est {\it tendue} si aucun disque
plong\'e dans $V$ ne s'appuie sur une courbe int\'egrale de $\xi$
transversalement \`a $\xi$. Dans le cas contraire, on dit que
$\xi$ est {\it vrill\'ee}.
Une structure de contact est {\it universellement tendue}
lorsque son rappel dans le rev\^etement universel de $V$ est tendu.

Soit  $\gamma$ une courbe plong\'ee dans $V$
qui borde une surface compacte, plong\'ee et orientable $S$.
Si $\gamma$ est transversale
\`a $\xi$, l'autoenlacement $l(\gamma )$ de
$\gamma$ est l'enlacement entre $\gamma$
et toute courbe obtenue en poussant l\'eg\`erement
$\gamma$ par une section non
singuli\`ere de $\xi \vert_S$. De m\^eme, si $\gamma$ est
{\it legendrienne}, c'est-\`a-dire tangente \`a $\xi$,
son invariant de Thurston-Bennequin est l'entier $tb(\gamma )$
obtenu en comptant l'enlacement entre $\gamma$ et sa d\'eformation dans
la direction d'un vecteur normal \`a $\xi$. Ces enlacements sont calcul\'es
avec l'orientation de $V$ qui rend $\xi$ positive.
Lorsque $\xi$ est tendue, si $\chi (S)$ d\'esigne la caract\'eristique
de $S$, l'autoenlacement $l(\gamma )$ et l'invariant de Thurston-Bennequin
$tb(\gamma )$  v\'erifient les in\'egalit\'es
de Bennequin~\cite{Be}~:
$$ l(\gamma )\leq -\chi (S) \; \; {\rm et} \; \; tb(\gamma )\leq -\chi (S).$$

Si $S$ est une surface close orient\'ee,
munie d'une forme d'aire $\omega$, et plong\'ee dans
une vari\'et\'e de contact
$(V,\xi )$, dont l'orientation est donn\'ee par $\xi$,
on appelle {\it feuilletage caract\'eristique} de $S$,
not\'e $\xi S$,
le feuilletage int\'egral du champ de vecteurs $Y$
donn\'e par $$i_Y \omega =\alpha \vert_{TS}$$ qui dirige $\xi \cap TS$.
Ses singularit\'es sont les points $x\in S$ o\`u $\xi (x)=T_x S$.
\begin{remarque}
 On insiste sur le fait que pour
d\'eterminer l'orientation du feuilletage caract\'eristique d'une
surface $S$, on utilise
une orientation de l'espace ambiant diff\'erente suivant que la structure
de contact $\xi$ consid\'er\'ee est positive ou n\'egative. La divergence d'une
singularit\'e du feuilletage caract\'eristique sera donc toujours positive
l\`a o\`u les orientations de la  structure $\xi$  et de $TS$ co\"\i ncident,
ce ind\'ependemment du signe de $\xi$.
\end{remarque}

\subsection{Paires}

Soit $V$ une vari\'et\'e de dimension trois orient\'ee.
On appelle {\it paire} de structures de contact
la donn\'ee d'une structure $\xi_+$ positive
et d'une structure $\xi_-$ n\'egative.

Dans la suite, on ne consid\`ere que des structures $\xi_+$ et
$\xi_-$ coorient\'ees (et donc orient\'ees).
On appelle {\it point de contact} entre $\xi_+$ et $\xi_-$
un point $x\in V$ o\`u $\xi_+ (x)=\xi_- (x)$.
Les points de contact sont de deux sortes~: positifs si
les coorientations de $\xi_+$ et $\xi_-$ co\"\i ncident
en $x$, et n\'egatifs si elles sont oppos\'ees.
Pour une paire $(\xi_+ ,\xi_- )$, on note
$\Delta_+$ le lieu des contacts positifs,
$\Delta_-$ le lieu des contacts n\'egatifs et $\Delta =\Delta_+ \cup \Delta_-$.

Soient $\alpha_+$ et $\alpha_-$ des \'equations de
$\xi_+$ et $\xi_-$, toujours suppos\'ees positives sur un
vecteur normal direct a, respectivement, $\xi_+$ et $\xi_-$.
Il existe un unique champ de vecteurs
$X$ inclus dans $\xi_-
$ et qui v\'erifie l'\'equation
$$i_X d\alpha_- \vert_{\xi_-} =\alpha_+ \vert_{\xi_-}.$$
Il est lisse,  nul le long de $\Delta$
et  dirige $\xi_+ \cap \xi_-$ sur $V\setminus \Delta $.
Sa divergence dans la direction de $\xi_-$ pour la forme $d\alpha_-$
est bien d\'efinie et non nulle aux points o\`u $\xi_+ =\xi_-$~:
$\mathcal{L}_X d\alpha_- \vert_{\xi_-} =d(i_X d\alpha_- )\vert_{\xi_-} =d\alpha_+ \vert_{\xi_+}$.
En un point de $\Delta_+$ (resp. $\Delta_-$),
la divergence de $X$ dans la direction de $\xi_\pm$ est
positive (resp. n\'egative).

On note $D^2 =\{ (x,y)\in \R^2,\;
x^2 +y^2 \leq 1\}$ le disque unit\'e de $\R^2$.
Une paire $(\xi_+ ,\xi_- )$, d\'efinie sur une vari\'et\'e 
close (compacte sans bord) $V$,  sera dite {\it normale} si
l'ensemble  $\Delta$ est un entrelacs lisse plong\'e
dans $V$, transversal \`a $\xi_\pm$
sauf en un nombre fini de points $(x_i )_{1\leq i\leq n}$, et qui v\'erifie
les propri\'et\'es suivantes.
\begin{itemize}
\item Si $a_i$ d\'esigne un arc ouvert transversal
\`a $\xi_\pm$ et  d\'elimit\'e
par $x_i$ et $x_{i+1}$ dans $\Delta$, alors  $a_i$
poss\`ede un voisinage tubulaire  $D^2 \times a_i$, $\{ (0,0)\} \times a_i =a_i$,
sur lequel
$\xi_+ \cap \xi_- \subset TD^2 \times \{ pt\}$ et
le feuilletage caract\'eristique  $\xi_+   D^2 \times \{ pt\}$ 
est \'egal au feuilletage
$\xi_-  D^2 \times \{ pt\}$ (\'egalit\'e dans $\Delta_+$, \'egalit\'e \`a orientation pr\`es
dans $\Delta_-$) et
est un feuilletage de $D^2 \times \{ pt\}$
par des selles ou par des foyers radiaux.

\item La sous-vari\'et\'e $\Delta$ a un contact quadratique avec $\xi_\pm$
en les $x_i$, et $x_i$ poss\`ede un voisinage $D^2 \times [-1,1]$,
$x_i =(0,0,0)$, o\`u $\xi_+ \cap \xi_- =\xi_+ \cap  TD^2 \times \{
t\} = \xi_- \cap TD^2 \times \{ t\}$, et o\`u le feuilletage
caract\'eristique de $\xi_+  D^2 \times \{ t\} =\xi_- D^2 \times
\{ t\} $, pour $t\in [-1,1]$, est le film d'une \'elimination
entre un foyer et une selle. En particulier, le feuilletage de
$D^2\times \{ 0\}$ poss\`ede une singularit\'e de type
naissance-mort en $(0,0,0)=x_i$. Dans ces coordonn\'ees, le champ
de vecteurs $X$ est tangent au champ de droites $\{ dt=dy=0\}$ sur
$\{ y=0 \} \times [-1,1]$
 (l'\'elimination se fait
le long de ces caract\'eristiques).
\end{itemize}

On montre que toute paire se laisse d\'eformer en une paire  normale.

\begin{proposition}\label{proposition : generique}
 Sur  une vari\'et\'e  close $V$, pour toute paire $(\xi_+ ,\xi_- )$
il existe une paire normale $(\xi'_+ ,\xi'_- )$ o\`u $\xi'_\pm$
est isotope \`a $\xi_\pm$.
\end{proposition}
\begin{proof}[D\'emonstration.]
On fixe une trivialisation de $TV$, ainsi qu'une m\'etrique sur $V$.
La donn\'ee de $\xi_+$ et $\xi_-$ est la donn\'ee de
deux applications $f_\pm : V\rightarrow S^2 \subset \R^3$.
La condition de contact est ouverte pour la topologie $C^1$, et
m\^eme, d'apr\`es le th\'eor\`eme de Gray, deux applications
$f$ et $g$ qui sont  $C^1$-proches donnent des structures de contact
isotopes.
On peut donc, quitte \`a perturber $\xi_+$ par isotopie,
se placer dans la situation g\'en\'erique (*) o\`u
l'application $F: V \rightarrow S^2 \times S^2$,
$F(x)=(f_+ (x),f_- (x))$, est transversale \`a
la diagonale et \`a l'antidiagonale de $S^2 \times S^2$.
Sous cette hypoth\`ese, $\Delta_+$ et $\Delta_-$
sont des sous-vari\'et\'es de dimension $1$ de $V$,
et le long de $\Delta$ la diff\'erentielle de $X$
est de rang $2$.
On suppose \'egalement v\'erifi\'ee la propri\'et\'e
g\'en\'erique~: les points
o\`u $\Delta_+$ et $\Delta_-$ sont tangents \`a $\xi_\pm$
sont des contacts quadratiques et donc isol\'es.
Les composantes de $\Delta \setminus (\cup_{1\leq i\leq n} x_i )$
sont de deux types suivant le signe du d\'eterminant
de la diff\'erentielle de $X$ dans la direction normale~:
{\it branches de foyers} si ce d\'eterminant est positif et
{\it branches de selles} s'il est n\'egatif. 

On suppose que $\Delta_+$ a des contacts quadratiques avec $\xi_\pm$~;
dans le cas contraire, l'\'etude se simplifie. Un  voisinage $N(\Delta_+ )$ de $\Delta_+$ est dit {\it normal}
s'il est constitu\'e des \'el\'ements suivants~:
\begin{itemize}
\item un voisinage $N(x_i )$ de chaque point de contact
quadratique $x_i$, $i=1,\dots ,n$, de la forme $N(x_i ) =\{ \vert
x\vert \leq z_0 , \vert y\vert \leq \sqrt{z_0} +z_0 ,\vert z\vert
\leq z_0\}$, o\`u $z_0
>0$ et $\xi_\pm (x_i )=\ker dz$. Le champ de vecteurs $X$
est transversal aux faces {\it verticales} $\{ \vert x\vert =z_0
\}$ et $\{ \vert y\vert =\sqrt{z_0} +z_0 \}$. Il 
est {\it horizontal} (i.e. inclus dans $\ker dz$) et sort de $N(x_i )$
le long de ces faces verticales,
sauf le long de la face $\{ y=\sqrt{z_0} +z_0 \}$ o\`u il rentre
dans $N(x_i )$. 
Pour tout $c\in [-z_0 ,z_0 ]$, 
l'arc $\{ \vert x\vert =z_0 , z=c\} \cup \{
y=-\sqrt{z_0} -z_0 ,z=c\}$ (orient\'e comme
bord du disque horizontal, lui-m\^eme coorient\'e par $\partial_z$, $\{ z=c\}$) 
est positivement transversal \`a
$\xi_+$ et, sauf aux deux points $\{ x=\pm z_0 ,y=0 ,z=-z_0 \}$,  
son int\'erieur est n\'egativement transversal \`a $\xi_-$.
 L'arc $\{ y=\sqrt{z_0} +z_0 ,z=c\}$ est
n\'egativement transversal \`a $\xi_+$ et tangent \`a $\xi_-$. En particulier,
du fait que $X$ y est horizontal,
$\xi_-$ a pour \'equation $dz=0$ le long de cet arc. 
La structure $\xi_-$ est \'egalement d'\'equation $dz=0$ le long de
l'arc $\{ y=0, z=-z_0 \}$, tandis que $\xi_+$ lui  est transversal.
Pour finir, le champ $X$ 
dirige le feuilletage caract\'eristique de la face
sup\'erieure $\{ z=z_0 \}$ pour $\xi_+$ et $\xi_-$ (seulement
en dehors de l'arc $\{ y=\sqrt{z_0} +z_0 \}$ pour cette
derni\`ere, qui est, d'apr\`es ce qui pr\'ec\`ede,  une ligne de singularit\'es de $\xi_- \{ z=z_0 \}$). 
\item un voisinage tubulaire $N =\cup_{1\leq j\leq n} N_j$ de 
la collection d'arcs $\Delta_+ \setminus \cup_{1\leq i\leq n} \Int (N(x_i ))$, 
o\`u les voisinages $N_j \simeq [-1,1]\times [-1,1] \times [0,m]$ 
sont deux \`a deux disjoints. Chaque $N_j$ est muni de
coordonn\'ees $(x,y,z)$ avec  $\Delta_+ \cap N_j=(0,0)\times [0,m]$. Il
intersecte $N(x_i)$, $i=1,\dots ,n$, le long du disque $\{ -z_0
-\sqrt{z_0} \leq y\leq 0 ,z=-z_0\} \subset N(x_i )$, du disque $\{
0\leq y\leq \sqrt{z_0} +z_0 \} \subset N(x_i )$, ou de l'ensemble
vide. S'il s'agit du voisinage d'une branche de foyers, il est
feuillet\'e par des disques horizontaux $\{ z=c\}$ dont les bords
sont positivement transversaux \`a $\xi_+$ et n\'egativement
transversaux \`a $\xi_-$ (le champ de vecteurs $X$ est 
horizontal et sortant le
long du bord), sauf en $\{ z=0 \}$ et $\{ z=m\}$ o\`u une ar\^ete
est tangente \`a $\xi_-$. S'il s'agit du voisinage d'une branche
de selles, le champ $X$ est horizontal et rentrant le long des faces verticales
$\{ \vert y\vert =1\}$, et horizontal et 
sortant le long des faces verticales $\{
\vert x\vert =1\}$. La structure $\xi_+$ est n\'egativement
transversale aux arcs $\{ \vert y\vert =1 ,z=c\}$, tandis que
$\xi_-$ leur est tangente, et positivement transversale aux arcs
$\{ \vert x\vert =1 ,z=c\}$, tandis que $\xi_-$ est
n\'egativement transversale \`a leur int\'erieur.
\end{itemize}
En renversant le sens de $X$, on d\'efinit de fa\c con similaire
une notion de voisinage normal pour $\Delta_-$.

\begin{lemme} Si $(\xi_+ ,\xi_- )$ est une paire de structures de contact
g\'en\'erique pour laquelle $\Delta$ est une sous-vari\'et\'e de dimension $1$
qui a  des contacts quadratiques avec $\xi_\pm$, alors  $\Delta$
poss\`ede un voisinage normal.
\end{lemme}
\begin{proof}[D\'emonstration.] On explique comment trouver un
voisinage normal de $\Delta_+$. L'\'etude pour $\Delta_-$ est
similaire.

Soit  $p$ un point de contact quadratique entre $\Delta_+$ et $\xi_\pm$.
On prend un voisinage $U_0 \subset \R^3$ de
$p=(0,0,0)$, muni de ses
coordonn\'ees cart\'esiennes $(x,y,z)$, pour lequel $\xi_\pm (p)=\{ dz=0\}$.
Par commodit\'e, on se r\'ef\'erera
\`a la coordonn\'ee $z$ comme \`a la direction {\it verticale}, et
aux coordonn\'ees $x$ et $y$ comme aux directions {\it
horizontales}.
Soit $X$ le champ de vecteurs lisse, qui dirige
$\xi_+ \cap \xi_-$ en dehors de $\Delta$ et nul sur $\Delta$,
d\'efini par l'\'equation $i_X d\alpha_- \vert_{\xi_-} =\alpha_+
\vert_{\xi_-}$. On \'ecrit le d\'eveloppement de $X$ \`a l'ordre
$1$  au voisinage de $(0,0,0)$~: $X=X_1 +o(\Vert (x,y,z)\Vert )$.
Comme $X$ est dans $\xi_+$, il v\'erifie une \'equation
$dz=fdx+gdy$, avec $f(0,0,z)=g(0,0,z)=0$. En particulier,
$X_1$ est horizontal~: $X_1 = (a_0 x+b_0 y +c_0 z)\partial_x
+(a_1 x+ b_1 y+c_1 z)\partial_y$.
Comme le rang de $DX_1$ en $p$ est $2$, on obtient,
apr\`es un changement de variables dans $T_p V$ qui donne une
forme r\'eduite de Jordan \`a $D_p X$, que $X_1$
est du type $X_1 =ax\partial_x +bz\partial_y$. Une fois choisi
$(\partial_x ,\partial_y ,\partial_z )$ dans $T_p V$, on \'etend
ce rep\`ere par de nouvelles coordonn\'ees $(x,y,z)$ de sorte que 
$\Delta_+ 
= \{ z=-y^2 ,x=0\}$ et $\xi_- =dz-xdy$ pr\`es de $p$.

Le d\'eveloppement de $X$ \`a un ordre sup\'erieur est, dans ces
coordonn\'ees~:
$$X =(x(R_1 (x,y,z ) +(z+y^2)R_2 (x,y,z)) \partial_x + ((z+y^2)R_3
(x,y,z)+ xR_4 (x,y,z)) \partial_y +$$$$ (bx(z+y^2 )+x^2 R_4
(x,y,z))\partial_z +o(\Vert (x,y,z)\Vert^2 )$$ o\`u~:
\begin{itemize}
\item $R_1$, $R_2$, $R_3$ et $R_4$ sont d'ordre $1$~;
\item $R_1$ et $R_3$ ont un terme constant non nul (resp. $a$ et $b$), \`a l'inverse de
$R_2$ et $R_4$~;
\item les termes d'ordre $3$ provenant de $y^2 R_2$ et $y^2 R_3$ sont
ceux apparaissant dans le d\'eveloppement de $X$ \`a l'ordre $3$.
\end{itemize}

 Pour fixer les id\'ees, on se place dans le cas
o\`u $b>0$, ce qui signifie que $\Delta_+$ se situe dans la zone
$z\leq 0$.

On consid\`ere alors le voisinage $U$  de $p$ donn\'e par $\{
\vert x\vert \leq  z_0  ,\vert y\vert \leq \sqrt{z_0} +z_0 ,\vert
z\vert \leq z_0 \}$, o\`u $z_0 >0$ est choisi assez petit.

Dans l'\'ecriture de la coordonn\'ee de $X$ sur $\partial_x$, les
termes en $z^2$, $yz$ et $y^3$ sont toujours au maximum de l'ordre de
$(z_0 )^{1/2 +1}$, c'est-\`a-dire n\'egligeables devant $z_0$. De
m\^eme, lorsque $\vert y\vert =\sqrt{z_0} +z_0$, $z+y^2$ est
toujours sup\'erieur \`a un $O(z_0 )$~; c'est donc le terme
dominant dans le coefficient de $\partial_y$. On v\'erifie ainsi
que  pour $z_0$ assez petit, $X$ est transversal aux faces
verticales $\vert x\vert =z_0$ et $\vert y \vert =\sqrt{z_0}
+z_0$. Il sort de $U$ le long de $\vert x\vert =z_0$ car $a>0$
(c'est la divergence de $X$ dans $\xi_- (p)$), sort de $U$ le long
de $y=-(\sqrt{z_0} +z_0 )$ et rentre dans $U$ le long de
$y=\sqrt{z_0} +z_0$.

Lorsque le voisinage $U$ est assez petit ($z_0$ est assez petit),
$R_3 (x,y,z)$  vaut environ $b$. De plus, comme la composante sur
$\partial_z$ est environ $bx(z+y^2 )$, toute orbite partant dans $U$ de
l'altitude $z_0$ rencontre une des faces verticales (ou plut\^ot
leur prolongement dans la direction des $z$) \`a une altitude
sup\'erieure \`a $z_0 - O(z_0 \sqrt{z_0} )$. Pour $z_0$ assez
petit, une telle orbite sort \`a une altitude positive proche de
$z_0$. On prend dans la face $\{ y=\sqrt{z_0} +z_0  \}$ l'arc
$A=\{ z=\frac{3}{4} z_0 \}$. Il est transversal (positivement) \`a
$\xi_+$ et tangent \`a $\xi_-$. On consid\`ere la surface
constitu\'ee de la r\'eunion des orbites de $X$ passant par $A$. 
De son intersection
avec $U$, on ne garde que la composante qui contient $A$. La
discussion effectu\'ee plus haut assure que celle-ci est compacte,
a son bord dans les faces verticales, est parall\`ele \`a la face
$\{ z=z_0 \}$ et est contenue dans la r\'egion $z>0$. De plus, son
feuilletage caract\'eristique pour $\xi_+$ et $\xi_-$ est dirig\'e
par $X$ (en dehors de l'arc $A$). Elle d\'ecoupe $U$ en deux
composantes. On nomme $U'$ celle qui contient $p$. Remplacer $U$
par $U'$ revient \`a remplacer la face horizontale $\{ z=z_0 \}$ par
une face dont le feuilletage caract\'eristique pour $\xi_+$ et
$\xi_-$ est dirig\'e par $X$.

On change la structure produit de  $U'$ pour obtenir les
propri\'et\'es de tangence et de transversalit\'e des cercles
horizontaux $\partial U' \cap \{ z=c\}$ par rapport 
\`a $\xi_+$ et $\xi_-$. Pour cela, on utilise le fait
que l'holonomie de $\xi_+$ autour du bord vertical  de $U'$
est positive, tandis qu'elle est n\'egative pour $\xi_-$. Ici
l'holonomie de $\xi_\pm$ est d\'efinie comme l'application de
premier retour du feuilletage caract\'eristique $\xi_\pm \partial
U'$ sur une ar\^ete verticale, par exemple $\{ x=-z_0 ,
y=\sqrt{z_0} +z_0 \}$. Elle est positive (resp. n\'egative ou
nulle) quand l'orbite issue d'un point $q$ revient \`a une
altitude $z$ sup\'erieure (resp. inf\'erieure ou \'egale) \`a
celle de $q$.

Si on suit le feuilletage caract\'eristique de $\xi_-$ le long de
la face verticale $\{ y=z_0 +\sqrt{z_0} \}$, puis le feuilletage
de $\xi_+$ le long des autres faces verticales, on obtient une
holonomie encore sup\'erieure \`a  celle obtenue en suivant
toujours le feuilletage caract\'eristique de $\xi_+$ (car $X$
rentre dans $U'$ le long de $\{ y=z_0 +\sqrt{z_0} \}$ et sort de
$U'$ le long des autres faces verticales). \`A l'inverse, si on
suit le feuilletage de $\xi_-$, on obtient une holonomie
n\'egative. Par le th\'eor\`eme des valeurs interm\'ediaires, on
est assur\'e d'obtenir une holonomie nulle, c'est-\`a-dire un
feuilletage par cercles, en suivant le feuilletage de $\xi_-$ sur
la face $\{ y=z_0 +\sqrt{z_0} \}$ et en suivant une certaine
trajectoire restant dans le c\^one d\'elimit\'e par $\xi_+$ et
$\xi_-$ le long des autres faces verticales. On change les coordonn\'ees
pour que ces cercles deviennent les bords des disques horizontaux.
En particulier, on change ainsi la face
horizontale inf\'erieure
de $U'$, en la faisant s'appuyer sur un des cercles d\'ecrit plus
haut. 
L'argument d\'evelopp\'e ci-dessus
permet \'egalement de s'assurer que, apr\`es repara\-m\'etrisation,
la structure $\xi_-$ soit d'\'equation $dz=0$ le long de l'arc $\{ y=0, z=-z_0 \}$.
Le voisinage obtenu est $N(x_i )$.

Soit maintenant $p$ un point de $\Delta_+$ o\`u $\Delta_+$ est transversal
\`a $\xi_\pm$. Un voisinage de $p$ s'envoie sur $\R^3$ muni de 
coordonn\'ees cart\'esiennes $(x,y,z)$ dans lesquelles 
$$p=(0,0,0), \; \; \Delta_+ =\{ x=y=0\}\; \; {\rm et} \; \;
\xi_\pm (0,0,0)= \{ dz=0\}.$$  On \'ecrit \`a nouveau 
le d\'eveloppement de $X$ \`a l'ordre
$1$  au voisinage de $(0,0,0)$.
Cette fois, on obtient que $X_1 = (a_0 x+b_0 y)\partial_x
+(a_1 x+b_1 y)\partial_y$. On est toujours dans la situation g\'en\'erique
o\`u la diff\'erentielle de $X$, et donc de $X_1$, est de rang $2$
en $p$, c'est-\`a-dire que $a_0 b_1 -a_1 b_0 \neq 0$. Suivant le
signe de ce d\'eterminant, on est en pr\'esence d'une branche de
foyers ou de selles. On trouve facilement un voisinage $N$ de
cette branche avec les propri\'et\'es recherch\'ees
(voir la discussion du cas des voisinages $N(x_i )$ ci-dessus
pour la mise en bonne position de $\xi_+$ et $\xi_-$
par rapport au bord des disques horizontaux).

Enfin, on peut v\'erifier qu'il est possible
d'ajuster  les choix de $U' =N(x_i )$ et
de $N$ afin qu'ils s'intersectent comme dans la d\'efinition de
voisinage normal. La raison en est \`a nouveau que dans 
les voisinages $U$, la coordonn\'ee de $X$
sur $\partial_z$ est major\'ee par le produit de $z_0$ par la coordonn\'ee
de $X$ sur $\partial_y$, ce qui permet de consid\'erer
le vecteur $X$ comme horizontal et d'interpoler entre $N(x_i )$ et $N$
dans le voisinage $U$. Les d\'etails sont laiss\'es au lecteur. 
\end{proof}

La preuve de la proposition~\ref{proposition : generique} se
conclut en appliquant le lemme~\ref{lemme:remplacement}
ci-dessous.
\end{proof}

\begin{lemme}\label{lemme:remplacement}
Si le lieu des contacts $\Delta$ d'une paire $(\xi_+ ,\xi_- )$
poss\`ede un voisinage normal, alors $(\xi_+ ,\xi_- )$ se laisse
d\'eformer en une paire normale.
\end{lemme}
\begin{proof}[D\'emonstration.] On consid\`ere un voisinage
normal $N(\Delta )$ de $\Delta$. Il s'agit de modifier $\xi_+$ et
$\xi_-$ dans $N(\Delta )$ pour obtenir une paire normale. Dans la
suite, on se concentre sur un voisinage de $\Delta_+$. La
construction au voisinage de $\Delta_-$ est similaire.

Le lemme d'\'elimination des singularit\'es de~\cite{Gi1}
permet de d\'eformer $\xi_+$ relativement \`a $V \setminus \Int (N(\Delta ))$
pour que~:
\begin{itemize} 
\item si $N_j$, $j=1,\dots n$, est le voisinage d'une branche de foyers (resp. de selles),
la structure $\xi_+$ trace un foyer radial (resp. une selle) sur chaque disque horizontal
de $N_j$~;
\item sur chaque voisinage $N(x_i )$, $i=1,\dots ,n$,  
la structure $\xi_+$ trace sur le feuilletage  par disques horizontaux de $N(x_i )$ 
le film d'une \'elimination, l'\'elimination
se faisant le long des caract\'eristiques $\{ x=0, z=c\}$, $c\in [-z_0 ,z_0]$.
\end{itemize}

La deuxi\`eme \'etape est de remplacer $\xi_-$ \`a l'int\'erieur
de $N(\Delta )$ par une structure de contact n\'egative $\xi'_-$ qui trace
sur l'int\'erieur des disques horizontaux de $N(\Delta )$ le m\^eme feuilletage
que $\xi_+$.

On commence par \'etendre $\xi_- \vert_{V\setminus \Int N(\Delta )}$ aux voisinages
$N_j$ qui sont des branches de selles.
Sur $N_j =[-1,1] \times [-1,1]\times [0,m]$, le feuilletage caract\'eristique
$\xi_+ \{ z=c\}$ du disque horizontal $[-1,1] \times [-1,1]\times \{ c\}$
est donn\'e par le noyau d'une $1$-forme $\beta_c$ avec 
$d\beta_c \vert_{[-1,1] \times [-1,1]\times \{ c\}} >0$, 
qu'on peut supposer ind\'ependante de $c$ sur un voisinage $U$ du bord vertical.
On cherche \`a \'etendre $\xi_-$ sur $N_j$ par le noyau d'une \'equation
du type $$\alpha'_- =dz-f\beta_z =0,$$ o\`u $f: N_j \rightarrow \R$ est strictement
positive sur $\Int (N_j )$ et v\'erifie les conditions impos\'ees par $\xi_-$ au bord.
La condition de contact est donn\'ee par~:
$$(df\wedge \beta_z +f^2 \beta_z \wedge \dot{\beta_z} +fd_h\beta_z )\wedge dz >0,$$
o\`u $\dot{\beta_z}$ est la d\'eriv\'ee de $\beta_z$ par rapport \`a $z$
et $d_h$ est la diff\'erentielle dans la direction horizontale.

Le point cl\'e est que, puisque $\xi_-$ est tangente aux
disques horizontaux $\{ z=c\}$ le long des  arcs
$\{ \vert y\vert =1 ,z=c\}$, la fonction $f$ est nulle
le long de ces m\^emes arcs. Elle peut donc \^etre choisie arbitrairement
petite et de diff\'erentielle  nulle en dehors de $U$. Comme sa valeur
est positive le long des ar\^etes $\{ \vert x\vert =1 ,z=c\}$, on peut \'egalement imposer
\`a sa d\'eriv\'ee dans la direction du feuilletage caract\'eristique de $\xi_+$ d'\^etre positive. Pour un tel
choix de fonction $f$, la forme $\alpha'_-$ est de contact~:
en dehors de $U$, $df$ est nulle et $f^2$
est n\'egligeable devant $f$, donc  
$(df\wedge \beta_z +f^2 \dot{\beta_z} \wedge \beta_z +fd_h\beta_z )\wedge dz$ est de l'ordre de
$fd_h\beta_z \wedge dz>0$. 
Dans $U$, $\dot{\beta_z }=0$ et $df\wedge \beta_z \wedge dz >0$,
tout comme $fd_h \beta_z \wedge dz >0$. L'in\'egalit\'e de
contact est donc encore satisfaite.

On s'occupe ensuite de la partie $N_f$ de $N$ qui est un voisinage
des branches de foyers. Une petite difficult\'e technique
appara\^\i t~: les bords des disques $[-1,1]\times [-1,1] \times
\{ 0\}$ et $[-1,1]\times [-1,1] \times \{ m\}$ ont une ar\^ete
legendrienne pour $\xi_-$. On commence donc par \'elargir $N_f$ en
$N_f'$ pour que $N_f'$ soit feuillet\'e par des disques dont les bords
sont transversaux positivement \`a $\xi_+$ et
n\'egativement \`a $\xi_-$ (et avec le bord vertical de
$N_f'$ transversal \`a $X$). On \'etend les coordonn\'ees de $N_f$ \`a 
$N_f'$ pour que
$\xi_+$ trace toujours un feuilletage radial sur chaque disque horizontal
de $N_f'$.
On remplace
alors  la structure $\xi_-$ sur $N_f'$ par une structure $\xi_-'$
tangente aux feuilletages caract\'eristiques radiaux trac\'es 
par $\xi_+$, en la faisant pivoter de sa
position initiale au bord jusqu'\`a la position horizontale, atteinte au
centre.

Reste \`a \'etendre $\xi_-$  dans les voisinages $N(x_i )$, $i
=1,\dots ,n$, des contacts quadratiques. Dans chacun des
voisinages $N(x_i )$, le feuilletage caract\'eristique trac\'e par
$\xi_+$ sur les disques horizontaux donne le film d'une \'elimination
d'un foyer avec une selle.
La remarque cl\'e est, comme dans le cas des selles, qu'on peut alors
\'etendre la restriction de $\xi_-$ au bord vertical de 
$N(x_i )$ en une structure de contact qui trace sur
l'int\'erieur de chacun des disques horizontaux le m\^eme feuilletage que
$\xi_+$, et ce car $\xi_-$ est
tangente aux arcs $\{  y =z_0 +\sqrt{z_0} ,z=c \}$, $c\in [-z_0 ,z_0]$. 
Ce proc\'ed\'e fournit une premi\`ere structure qui ne co\"\i ncide pas 
n\'ecessairement avec $\xi_-$
le long des faces horizontales $\{ z=\pm z_0 \}$. Pour obtenir une structure
qui \'etend $\xi_-$, on interpole entre cette premi\`ere extension, d'\'equation
$dz+\beta_z =0$, et $\xi_-$, d'\'equation $dz+\delta_z =0$ ($\beta_z$ et $\delta_z$
ont m\^eme noyau dans $\Int (N(x_i ))$),
dans un voisinage  de $\{ z=\pm z_0 \}$. On  obtient  un champ de plans 
avec une \'equation du type $dz +\chi (z)\beta_z +(1-\chi (z))\delta_z =0$,
o\`u $\chi$ est une fonction de coupure,
qui est une extension de $\xi_-$ et qui va automatiquement \^etre de contact.

Pour finir, on v\'erifie que la structure  $\xi_-'$
produite  est isotope \`a la structure initiale 
$\xi_-$. Pour cela, on observe par exemple que les structures
$\xi_-$ et $\xi_-'$ sont tendues sur $N(\Delta )$ (car
\guil{horizontales}). La structure $\xi_-'$ est \'egale \`a
$\xi_-$ le long de $\partial N(\Delta )$ et toutes deux
poss\`edent un cercle m\'eridien transversal et d'autoenlacement
$-1$. L'isotopie est alors donn\'ee par un avatar du th\'eor\`eme
de classification d'Eliashberg~\cite{El} (il est \'egalement possible ici
de trouver une isotopie explicite par des arguments \'el\'ementaires).
\end{proof}

\begin{remarque} La perturbation d'une paire g\'en\'erique $(\xi_+ ,\xi_- )$
en une paire normale $(\xi_+' ,\xi_-' )$ peut \^etre r\'ealis\'ee
en fixant $\Delta$.
\end{remarque}

Soit $(\xi_+ ,\xi_- )$ une paire de structures de
contact normale. On note toujours $X$ un
champ de vecteurs lisse, nul sur $\Delta$
et qui dirige le champ de droites orient\'ees $\xi_+ \cap \xi_-$
sur $V\setminus \Delta $.
Comme dans \cite{ET}, on peut alors d\'efinir
quatre champs de plans $\lambda_\pm^{\pm \infty}$
en posant, pour $x \in V\setminus \Delta $~:
$$\lambda_\pm^{\pm \infty} (x) =\lim_{t\rightarrow \pm \infty}
X^t_* (X^{-t} (x)) \xi_\pm (X^{-t} (x)).$$
Ici, $X^t$ d\'esigne le flot de $X$ \`a l'instant $t$.
On pose de plus, pour $x\in \Delta_+$, $$\lambda_\pm^{+\infty} (x)=\xi_+ (x)=
\xi_- (x),$$ et pour $x\in \Delta_-$, $$\lambda_\pm^{-\infty} (x)=\xi_+ (x)=
\xi_- (x).$$

\begin{proposition}\label{proposition : egalite}
Si $V$ est close, si $(\xi_+ ,\xi_- )$ est une paire normale 
de structures de contact sur $V$ et
si $\Delta_-$ est transversale \`a $\xi_+$ (s'il n'y a
pas de contact quadratique le long de $\Delta_-$), alors
$\lambda_+^{+\infty} =\lambda_-^{+\infty} =:
\lambda^{+\infty}$ sur $V\setminus \Delta_-$.
\end{proposition}
\begin{proof}[D\'emonstration.]
On reprend les notations de \cite{ET}.
On fixe une m\'etrique $g$ sur $V$ et
on note $\nu$ le  champ de plans orthogonal \`a $X$
sur $V\setminus \Delta $.
On pose $\lambda_\pm =\xi_\pm \cap \nu$, et,
pour $t\in \R$ et $x \in V\setminus \Delta $,
$$\lambda^t_\pm (x) = X^t_* (X^{-t} (x)) \xi_\pm (X^{-t} (x))\cap \nu (x).$$

On notera de la m\^eme fa\c con la droite $\lambda^t_\pm (x)$
et le plan contenant cette droite et le vecteur $X$.
Si $\theta^t_\pm (x)$ d\'esigne l'angle (calcul\'e avec $g$)
entre $\lambda_\pm^t (x)$ et $\lambda_\pm (x)$,
les conditions de contact imposent~:
$$\frac{d\theta_+^t (x)}{dt} >0 \; {\rm et} \; \frac{d\theta_-^t (x)}{dt} <0.$$
Ces conditions, plus le fait que les fonctions
$\theta_\pm^t (x)$ soient born\'ees
(\`a $x$ fix\'e, et ind\'ependem\-ment de $t$), suffisent pour assurer
l'existence de limites quand $t$ tend vers $\pm \infty$.
Il est alors automatique que les plans  limites $\lambda_\pm^{\pm \infty}$
sont invariants par le flot de $X$.

On traite le cas o\`u $t$ tend vers $+\infty$.
On note $O$ l'orbite de $X$ qui passe par un point
$x \in V\setminus \Delta $.
Pour prouver la proposition~\ref{proposition : egalite} dans cette situation,
on distingue six  cas. Les deux premiers sont essentiellement trait\'es dans
\cite{ET}. Pour le confort du lecteur, on reproduit tout de m\^eme
la preuve de ces deux cas sp\'eciaux. La preuve propos\'ee pour le
cas 2 tient compte du fait qu'ici, contrairement \`a~\cite{ET},
l'espace $V\setminus \Delta $ n'est pas compact.
Comme la paire consid\'er\'ee est normale, soit l'orbite $O$ a
une limite $x_{-\infty}$ dans $\Delta$ lorsque $t$ tend vers $-\infty$,
soit elle a un point d'accumulation en temps n\'egatif
dans $V\setminus \Delta$.\\
\\
{\bf Cas 1.} $O$ est une orbite p\'eriodique.

On note $L:\nu \rightarrow \nu$ la diff\'erentielle de l'application de
premier retour du flot de $X$ en $x$. Les plans $\lambda_\pm^{+\infty} (x)$
et $\lambda_\pm^{-\infty} (x)$ sont distincts et invariants
par le flot de $X$. Si $\lambda_-^{+\infty} \neq \lambda_+^{+\infty}$,
l'application lin\'eaire $L$ a au moins trois directions propres distinctes
et se trouve donc \^etre une homoth\'etie. Ceci contredit le fait
que les directions $\lambda_\pm$
de $\xi_\pm $ ne sont pas laiss\'ees invariantes
par $L$.
On peut \'egalement argumenter de la fa\c con suivante, ce qui sera
utile pour l'\'etude du cas 2.
Les directions $\lambda_\pm^{+\infty} (x)$ sont des directions propres de $L$ et
les droites $\lambda_\pm (x)$ sont situ\'ees dans un m\^eme
secteur de $\nu \setminus (\lambda_+^{+\infty} (x) \cup \lambda_-^{+\infty} (x) )$.
Comme telles, ces derni\`eres doivent \^etre tourn\'ees dans le m\^eme
sens par l'application lin\'eaire $L$, ce qui n'est pas le
cas, puisque
$\frac{d\theta_+^t (x)}{dt} >0 \; {\rm et} \; \frac{d\theta_-^t (x)}{dt} <0.$
$ $\\
\\
{\bf Cas 2.} $O$ a
une valeur d'adh\'erence en temps n\'egatif dans $V\setminus \Delta $.

Dans ce cas, on fait fonctionner une version quantitative du cas 1).
On suppose que $\lambda_+^{+\infty} (x)\neq \lambda_-^{+\infty} (x)$.
Soit $x_\infty \in V\setminus \Delta $
la limite d'une suite $x_n =X^{t_n} (x)$ de points de $O$, avec $t_n \rightarrow -\infty$.
Il existe $\alpha >0$ tel que, si $t$ est assez grand, $\theta_+^t (x_\infty ) \geq \alpha$ et
$\theta_-^t (x_\infty )\leq -\alpha$.
On en d\'eduit qu'il
existe $\alpha >0$, un voisinage $U$ de $x_\infty$ et $T>0$, tels
que pour tout $p\in U$ et pour tout $t\geq T$,
$\vert \theta_\pm ^t (p)\vert >\alpha$. Quitte \`a prendre $n$ assez
grand, on suppose que tous les $x_n$ sont dans $U$. On identifie
$U$ \`a $\R^3$, en rendant $\nu$ tangent \`a $T(\{ *\} \times \R^2 )$
et $X$ colin\'eaire \`a $\R \times (*,*)$, si bien qu'on identifie
tous les plans de $\nu$ par translation.

On suppose pour l'instant qu'on peut trouver une suite $(x_{\phi (n)} )$ extraite de
$x_n$ telle que $\lambda_+^{+\infty} (x_{\phi (n)} )$ et $\lambda_-^{+\infty} (x_{\phi (n)} )$
aient pour limites respectives $\mu_+$ et $\mu_-$, lorsque $n$ tend vers $+\infty$, avec
$\mu_+ \neq \mu_-$.
Si $n$ et $p$ sont assez grands avec $n<p$, il existe dans ce cas
une application lin\'eaire
$A_{n,p} :\R^2 \rightarrow \R^2$ qui est $\epsilon$-proche de l'identit\'e
($\epsilon \ll \alpha$) et qui envoie
$\lambda_\pm^{+\infty} (x_{\phi (p)} )$ sur $\lambda_\pm^{+\infty} (x_{\phi (n)} )$.
Mais alors, d\`es que $t_{\phi (n)} -t_{\phi (p)} >T$,
la compos\'ee de $A$ avec la diff\'erentielle
$X_*^{t_{\phi (p)} -t_{\phi (n)}} (x_{\phi (n)} )$ a $\lambda_\pm^{+\infty} (x_{\phi (n)} )$
comme espaces propres, et  fait tourner  les deux
droites $\lambda_{\pm} (x_{\phi (n) } )$,
situ\'ees dans la m\^eme composante de
$\R^2 \setminus (\lambda_+^{+\infty} (x_{\phi (n)} )
\cup \lambda_-^{+\infty} (x_{\phi (n)} ))$, dans des sens oppos\'es
(car $\epsilon \ll \alpha$), ce qui est impossible.

Il reste \`a \'etudier le cas o\`u $\mu_+ =\mu_-$.
Cette fois, on peut trouver $n<p$ pour que les points $x_{\phi (n)}$ et
$x_{\phi (p)}$ soient dans $U$ et que l'angle form\'e par 
les droites  $\lambda_+^{+\infty}(x_{\phi (p)} )$ et
$\lambda_-^{+\infty}(x_{\phi (p)} )$ soit tr\`es petit (quitte \`a prendre
$p$ assez grand) devant l'angle
form\'e par $\lambda_+^{+\infty}(x_{\phi (n)} )$
et $\lambda_-^{+\infty}(x_{\phi (n)} )$. 
On peut alors trouver une application  lin\'eaire
(pas proche de l'identit\'e) $A_{n,p}$
qui envoie les deux droites $\lambda_\pm^{+\infty}(x_{\phi (p)} )$ sur
$\lambda_\pm^{+\infty}(x_{\phi (n)} )$ et qui \guil{\'ecarte} les
droites  $\lambda_\pm (x_{\phi (p)} )$~: les droites $\lambda_\pm (x_{\phi (n)} )$
sont dans un m\^eme secteur de 
$$\R^2 \setminus (A_{n,p} \lambda_+ (x_{\phi (p)} )
\cup A_{n,p} \lambda_- (x_{\phi (p)} )).$$
Par les conditions de contact, le secteur
angulaire d\'elimit\'e par les droites $\lambda_+ (x_{\phi (n)} )$
et $\lambda_- (x_{\phi (n)} )$ est envoy\'e par la diff\'erentielle
$X_*^{t_{\phi (p)} -t_{\phi (n)}} (x_{\phi (n)} )$ sur un secteur contenant
les droites $\lambda_\pm (x_{\phi (p)} )$.
La composition de $A_{n,p}$ avec la diff\'erentielle
$X_*^{t_{\phi (p)} -t_{\phi (n)}} (x_{\phi (n)} )$
a $\lambda_\pm^{+\infty} (x_{\phi (n)} )$ comme espaces propres et 
fait comme pr\'ec\'edemment tourner les
droites $\lambda_\pm (x_{\phi (n)} )$ en sens oppos\'es. C'est la contradiction
recherch\'ee.
$ $\\
\\
{\bf Cas 3.} $O$  provient
d'un foyer de $\Delta_+$.

Dans ce cas, on peut mod\'eliser un voisinage
du foyer $e$ en question par $$\R^3 =\{ (r,\theta ,z)\in ]0,\epsilon ]\times
\R /\Z \times [-1,1]\},$$
$\xi_+ =\{ dz-f(r^2 ,\theta ,z)d\theta =0 \}$, $f(0,\theta ,z)=0$
et $\xi_- = \{ dz-g(r^2 ,\theta ,z)d\theta =0 \}$, $g(0,\theta ,z)=0$.
On v\'erifie alors que dans ces coordonn\'ees,
$\lambda_\pm^{+\infty} ( r,\theta ,z)=\{ dz=0\}$.
Cette \'egalit\'e, obtenue explicitement au voisinage de $e$,
se propage en $x$ par le flot de $X$ (de m\^eme que
la diff\'erentiabilit\'e de ce champ de plans limite).

$ $\\
{\bf Cas 4.} $O$ est la s\'eparatrice instable d'une selle de
$\gamma_+$.

On utilise un mod\`ele comme dans le cas $3$, qui montre
que le long de la s\'eparatrice instable, au voisinage
de la selle $h$, les plans limites $\lambda_\pm^{+\infty}$
co\"\i ncident avec le plan tangent au feuilletage
par disques qui contient $X$.
Pr\'ecisemment, si on conjugue un voisinage de $h$
\`a sa lin\'earisation, on obtient, dans des coordonn\'ees
$(a,b,c)\in D^2 \times [-\epsilon ,\epsilon ]$, que
le flot de $X$ est $X^t (a,b,c) =(a\exp (\mu_+ t),b\exp (\mu_- t),c)$.
Dans ce mod\`ele, $\mu_+ >-\mu_- >0$ car $h$ est \`a divergence
positive et $x$ est sur la s\'eparatrice $s$
qui est $\{ a>0,b=c=0\}$ pr\`es de $h$. Les plans de contact $\xi_+$ et
$\xi_-$
ont quant \`a eux pour \'equations le long de $s$ respectivement
$dc+\alpha adb =0$ et $dc-\beta adb=0$.
On calcule alors que $X^t_* (X^{-t} (a,0,0)) \xi_+ (X^{-t} (a,0,0))$
est d'\'equation $dc+\alpha \exp ((-\mu_+ -\mu_- )t)db=0$.
Lorsque $t$ tend vers l'infini la limite est bien $\{ dc=0\}$.
On aboutit au m\^eme r\'esultat pour $\xi_-$.

Comme dans le cas pr\'ec\'edent, l'\'egalit\'e
$\lambda_+^{+\infty} =\lambda_-^{+\infty}$ se prolonge
en $x$ par le flot de $X$.
$ $\\
\\
{\bf Cas 5.} $O$ est la s\'eparatrice instable d'une selle
de $\Delta_-$.

Dans ce cas, $\lambda_\pm^{+\infty}$ est le plan tangent
\`a la r\'eunion des vari\'et\'es instables issues de $\Delta_-$. On d\'emontre
ce fait en reprenant le mod\`ele local du cas 4), avec cette
fois $-\mu_- >\mu_+ >0$.
$ $\\
\\
{\bf Cas 6.} $O$ provient d'une naissance-mort de $\Delta_+$.

Si $O$ provient de $\Delta_+$, et de l'int\'erieur de la partie
foyer de la naissance-mort,  le plan $\lambda_\pm^{+\infty}$
est \`a nouveau, au voisinage de la singularit\'e, le plan
tangent au feuilletage en disques, donn\'e par la forme
normale, qui contient $X$. Un germe de  feuille int\'egrale passant par $x$
est  en particulier lisse, comme dans le cas des foyers.
Si $O$ est l'une des deux orbites bord de la partie foyer,
l'\'ecriture d'un mod\`ele local fournit \'egalement l'\'egalit\'e recherch\'ee.
\end{proof}

\begin{proposition} Si $(\xi_+ ,\xi_- )$ est une paire normale
sur une vari\'et\'e close $V$  et si
$\Delta_-$ est transversal \`a $\xi_+$,
alors le  champ de plans $\lambda^{+\infty}$ est
continu sur $V\setminus \Delta_-$,  invariant par le
flot de $X$ et localement int\'egrable.
Lorsque de plus  $\Delta_+$ est transversal \`a
$\xi_+$, alors les deux champs de plans
$\lambda^{+\infty}$ et $\lambda^{-\infty}$
(d\'efini de fa\c con similaire \`a $\lambda^{+\infty}$) sont
transversaux en dehors de $\Delta$.
\end{proposition}
En particulier par tout point o\`u $\lambda^{\pm \infty}$ est
d\'efini passe un germe de  vari\'et\'e int\'egrale.
Ce germe n'est pas forc\'ement unique, et $\lambda^{\pm \infty}$
ne poss\`ede pas forc\'ement de feuilletage int\'egral global.

\begin{proof}[D\'emonstration.]
On montre que $\lambda^{+\infty}$ est continu sur $V\setminus
\Delta_- $. On note $\theta_+^{'t} (x)= \theta_+^t (x)$,
$\theta_-^{'t} (x)$ l'angle entre $\lambda_-^t (x)$ et $\lambda_+ (x)$, et 
$\theta^{'\infty}$ la limite commune
de $\theta_+^{'t}$ et $\theta_-^{'t}$ quand $t$ tend vers $+\infty$.
Comme $\theta_-^t (x)$, 
la fonction $\theta_-^{'t} (x)$ est une fonction d\'ecroissante de $t$.

Soit $\epsilon >0$ et $x\in V\setminus \Delta $.
Il existe $t_0$ tel que, pour $t\geq t_0$, $$\theta^{'\infty} (x) -\epsilon
<\theta_-^{'t} (x) <\theta_+^{'t} (x) <\theta^{'\infty} (x) +\epsilon$$
(pour des rel\`evements convenablement choisis des angles
$\theta'$ dans $\R$).
Par continuit\'e du flot de $X$ par rapport \`a $x$,
si $y$ est assez proche de $x$, on a
$\vert \theta_\pm^{'t_0} (y) -\theta_\pm^{'t_0} (x) \vert <\epsilon /2$.
Mais alors, pour $t\geq t_0$,
$\theta^{'\infty} (x) -3\epsilon/2 <\theta_-^{'t_0} (y) <
\theta_-^{'t} (y) <\theta^{'\infty} (y) <\theta_+^{'t} (y) <\theta_+^{'t_0} (y) <
\theta^{'\infty} (x) +3\epsilon /2$.
On en d\'eduit la continuit\'e 
de $\theta^{'\infty}$, et donc de $\lambda^{+\infty}$, sur $V\setminus \Delta$.
 
On a \'egalement bien \'evidemment continuit\'e aux points
de $\Delta_+$, car le plan $\lambda^{+\infty}$ est coinc\'e entre
$\xi_+$ et $\xi_-$ qui tendent tous deux vers le m\^eme
plan. C'est donc aussi le cas de
$\lambda^{+\infty}$.

On montre maintenant que $\lambda^{+\infty}$ est localement int\'egrable~:
par tout point de $V\setminus \Delta_-$ passe un germe de
surface int\'egrale.

Si $x \in V\setminus \Delta$, on se donne
un germe de disque $D$ centr\'e en  $x$ et transversal \`a $X$.
Le champ de plans $\lambda^{+\infty}$ intersecte le plan tangent \`a $D$
en un champ de droites continu. Par le th\'eor\`eme
de Cauchy-P\'eano, ce champ de droites admet
un germe de courbe int\'egrale $\gamma$ passant par $x$.
Le germe de surface obtenu en poussant $\gamma$ par le
flot de $X$ est tangent \`a $\lambda^{+\infty}$, puisque $X$ est
tangent \`a $\lambda^{+\infty}$ et que $\lambda^{+\infty}$ est invariant par le
flot de $X$.

Si $x$ est un foyer de $\Delta_+$, le mod\`ele local donn\'e
dans la preuve de la proposition~\ref{proposition : egalite}
donne que $\lambda^{+\infty}$ est lisse et uniquement int\'egrable au
voisinage de $x$ (c'est le champ de plans horizontal).

Si $x$ est une selle de $\Delta_+$, on se replace dans les coordonn\'ees
de la proposition~\ref{proposition : egalite}, o\`u le flot de
$X$ au temps $t$ dans un voisinage $U$ de $x$
est $(a,b,c)\rightarrow (a\exp (\mu_+ t),b\exp (\mu_- t),c)$,
et $x=(0,0,0)$.
On se donne deux germes de disques $D_1 \subset U$ et $D_2 \subset U$ transversaux
\`a $X$ et intersectant les deux s\'eparatrices stables de $x$.
On choisit  ensuite \`a l'aide du th\'eor\`eme de
Cauchy-P\'eano, pour $i=1,2$,
des courbes int\'egrales $\gamma_i \subset U$ de $TD_i \cap \lambda^{+\infty}$
passant par ces s\'eparatrices.
On consid\`ere alors l'ensemble $S$ des trajectoires de $X\vert_U$
issues de $\gamma_1$ et $\gamma_2$, augment\'e des s\'eparatrices
instables de $x$ dans $U$.
On montre que $S$ donne
 un germe de surface int\'egrale de
$\lambda^{+\infty}$ pr\`es de $X$.
Tout d'abord,  $S$ est le graphe d'une fonction $f$
continue~: $S=\{ (a,b,f(a,b))\}$ pr\`es de $x$.
Il est automatique que les ensembles $S_1$ et $S_2$ des trajectoires issues
respectivement de $\gamma_1$ et $\gamma_2$ forment deux surfaces int\'egrales ouvertes
de $\lambda^{+\infty}$. Il faut voir que ces deux surfaces se recollent
le long des s\'eparatrices instables de $x$ pour former
une surface int\'egrale globale. Pour cela, on utilise
le mod\`ele local~: une orbite issue d'un point $(a,b,\epsilon )$
passe en $(a', b(\frac{a}{a'} )^{-\frac{\mu_-}{\mu_+}}, \epsilon )$,
avec  $\frac{b}{a}(\frac{a}{a'} )^{-\frac{\mu_-}{\mu_+}} \rightarrow +\infty$
quand $a\rightarrow 0$, car $0<-\frac{\mu_-}{\mu_+} <1$.
Ceci permet de majorer le taux de variation de $f$ en $(a',0)$
par celui en un point $(0,b)$ o\`u $f$ est diff\'erentiable et
de diff\'erentielle nulle. On obtient que $f$ est diff\'erentiable
de diff\'erentielle nulle en $(a',0)$, ce qui confirme que
$S$ poss\`ede un plan tangent \'egal \`a $\lambda^{+\infty}$
le long des s\'eparatrices instables de $x$.

On traite de mani\`ere similaire le cas o\`u $x$ est une naissance-mort~:
on a d\'ej\`a une demi-surface int\'egrale du c\^ot\'e foyer,
qui se recolle avec une surface int\'egrale ouverte satur\'ee par $X$
du c\^ot\'e selle, comme dans le sous-cas pr\'ec\'edent (m\^eme
si cette fois le mod\`ele local ne donne pas une convergence vers $+\infty$,
mais de fa\c con suffisante un quotient de la deuxi\`eme
coordonn\'ee par $a$ minor\'e par une constante $>0$).
\end{proof}

\section{Paires positives de structures de  contact}

On suppose maintenant que $\Delta_- =\emptyset$. La paire de
contact $(\xi_+ ,\xi_- )$ est alors dite {\it positive}. On
rappelle qu'on ne consid\`ere que des paires normales. Pour
celles-ci, on s'int\'eresse au champ de plans $\lambda^{+\infty}$
qui est continu et d\'efini sur tout $V$. On le note $\lambda$
pour simplifier les notations. On poss\`ede de nombreux exemples
de telles paires de structures de contact positives. Un cas
important se trouve dans~\cite{ET}~: tout feuilletage $\F$,
diff\'erent du feuilletage en sph\`eres de $S^1 \times S^2$, peut
\^etre approch\'e par des structures de contact positives et
n\'egatives. Un tel couple de structures de contact proches de
$\F$ est automatiquement positif et se laisse d\'eformer
en une paire normale et positive d'apr\`es la proposition~\ref{proposition :
generique}.

On rappelle qu'un feuilletage $\F$ de codimension $1$ sur une
vari\'et\'e de dimension trois est {\it tendu} s'il poss\`ede
une transversale qui coupe toutes ses feuilles. L'existence
d'une telle transversale implique l'absence de {\it composante de Reeb}.

\begin{lemme}\label{lemme : uniforme}
Si $(\xi_+ ,\xi_- )$ est une paire normale et positive 
sur une vari\'et\'e close $V$, 
les champs de plans $\xi_\pm^t =X^t_* \xi_\pm$ convergent uniform\'ement vers
$\lambda$ lorsque $t$ tend vers $+\infty$.
\end{lemme}
\begin{proof}[D\'emonstration.] Les fonctions $\theta_\pm^t (x)$
sont monotones en $t$ et convergent simplement vers une limite
continue. Comme la vari\'et\'e $V$ est compacte, cette convergence est uniforme
d'apr\`es le th\'eor\`eme de Dini.
\end{proof}

\begin{proposition} Si $(\xi_+ ,\xi_- )$ est une paire normale positive et
si $\xi_+$ ou $\xi_-$ est tendue, alors $\lambda$ ne poss\`ede pas
de sph\`ere int\'egrale.
\end{proposition}
\begin{proof}[D\'emonstration.] Soit $e(\xi_\pm ) \in H^2 (V,\Z )$
la classe d'Euler de $\xi_\pm$. Si $S$ est une sph\`ere
int\'egrale de $\lambda$, alors le feuilletage caract\'eristique de
$S$ pour $\xi_+$ et $\xi_-$  ne poss\`ede que des singularit\'es positives.
On en d\'eduit que $e(\xi_\pm ).[S] >0$.
Or pour une structure $\xi$ tendue, on a $e(\xi ).[S]=0$
d'apr\`es \cite{Be, El}.
\end{proof}

\begin{remarque} Si $S$ est une surface int\'egrale de $\lambda$,
m\^eme si $S$ n'est pas lisse, son feuilletage caract\'eristique
est bien d\'efini puisque dirig\'e par le champ de vecteurs lisse $X$.
\end{remarque}

\begin{theoreme}\label{theoreme : reeb}
Si $(\xi_+ ,\xi_- )$ est une paire normale positive sur une vari\'et\'e close $V$, si
 $\xi_+$ et $\xi_-$ sont tendues  et si
$\lambda$ est uniquement int\'egrable, alors le feuilletage
int\'egral $\F$  de $\lambda$ ne poss\`ede pas de composante
de Reeb d'\^ame homologue \`a z\'ero dans $H_1 (V,\Q )$.
De plus, $V$ porte un feuilletage sans composante de Reeb et
en particulier son rev\^etement universel est $\R^3$.
\end{theoreme}
\begin{remarque} R\'eciproquement, si $\lambda$ a pour feuilletage
int\'egral un feuilletage tendu  $\F$,
les structures $\xi_\pm$ sont tendues (et m\^eme {\it symplectiquement
remplissables}),
puisque isotopes \`a des structures proches de $\F$
(lemme~\ref{lemme : uniforme}) et donc soumises aux conclusions
de~\cite{ET}.
\end{remarque}

\begin{proof}[D\'emonstration.]
Soit $\F$ le feuilletage int\'egral de $\lambda$.
On suppose que $\F$ contient une composante de Reeb,
c'est-\`a-dire en particulier une feuille torique $T$ qui borde un tore solide.

On se donne un m\'eridien $m$ trac\'e sur $T$ et
 une bande $B$  d'\^ame $m$, transversale \`a $\F$.
La remarque cl\'e est que
si l'holonomie de $\F$ le long de $m$ est non triviale
(si son germe n'est pas l'identit\'e),
on peut tracer sur $B$ une courbe plong\'ee $\gamma$
positivement transversale \`a $\F$ et isotope \`a $m$ dans $B$.
Le feuilletage trac\'e par $\F$ sur $B$ n'\'etant que $C^0$,
pour trouver cette transversale, on utilise ici de mani\`ere
cruciale que le champ de droites continu $T\F \cap TB$ est
uniquement int\'egrable, conjugu\'e \`a un r\'esultat
classique de Kneser (voir par exemple le \guil{th\'eor\`eme de
Kneser} dans  le livre~\cite{Har}).
Cette courbe $\gamma$ est le bord d'un disque plong\'e
dans $V$.

 Comme $\xi_\pm^t$ converge uniform\'ement vers
$\lambda$ lorsque $t$ tend vers $+\infty$,
il existe $t_0$ tel que si $t\geq t_0$,
les structures $\xi^t_\pm$ soient transversales
\`a $\gamma$, positivement pour les deux.
On en d\'eduit une homotopie entre $\xi_+^{t_0}$
et $\xi_-^{t_0}$ parmi les champs de plans transversaux
\`a $\gamma$. Elle est construite en concat\'enant le
chemin entre $\xi_+^{t_0}$ et $\lambda$ donn\'e par  les $\xi_+^t$,
$t\geq t_0$,
et le chemin entre $\lambda$ et $\xi_-^{t_0}$ donn\'e par les $\xi_-^t$,
$t\geq t_0$.

Cette homotopie indique que les autoenlacements de $\gamma$
pour les structures $\xi_+^{t_0}$ et $\xi_-^{t_0}$
sont oppos\'es (les structures sont de signes oppos\'es).
On ne peut donc pas avoir pour les deux structures
l'in\'egalit\'e de Bennequin~\cite{Be}~: $l(\gamma )\leq -1$.
Une au moins des structures est vrill\'ee.
Par contradiction, l'holonomie de $\F$  le long de $m$ est l'identit\'e.

Soit \`a pr\'esent $l$ une longitude trac\'ee sur $T$.
On montre que si $[l]$ est nulle dans $H_1 (V,\Q )$, alors
le germe de l'holonomie de $\F$ le long de $l$
est l'identit\'e au moins d'un c\^ot\'e de $T$.

Si tel n'est pas le cas, m\^eme lorsque $l$ n'est pas homologue \`a
z\'ero, on peut trouver un petit
voisinage tubulaire $T\times [-1,1]$ de $T$, $T=T\times \{ 0\}$,
o\`u les tores $T\times \{\pm 1\}$ sont transversaux \`a $\F$.
Comme l'holonomie de   $\F$ le long de $m$ est l'identit\'e,
on a de plus que les feuilletages trac\'es par $\F$
sur $T\times \{ \pm 1\}$ sont des feuilletages en cercles m\'eridiens.

\`A ce stade, on peut utiliser le fait, d\^u \`a Haefliger~\cite{Ha},
que l'union des feuilles compactes de $\F$ est un compact de $V$
pour consid\'erer
les feuilles toriques maximales parmi celles qui bordent un tore solide
(la maximalit\'e se rapporte \`a l'inclusion des tores solides).
D'apr\`es ce qui pr\'ec\`ede, un tel tore maximal est inclus dans un tore
solide sur le bord duquel $\F$ trace un feuilletage  par cercles m\'eridiens.
On peut remplacer $\F$ sur ce tore par un feuilletage en disques.
Cette construction conduit \`a un feuilletage $\F'$ sans composante de Reeb.

Pour revenir \`a l'\'etude du feuilletage
$\F$ dans le cas o\`u $[l]=0$, on note que,
comme $\xi_+^t$ converge uniform\'ement vers $\F$,
pour $t$ assez grand, $\xi_+^t$ est transversal \`a $T\times \{ \pm1\}$
et le feuilletage caract\'eristique $\xi_+^t T\times \{ \pm 1\}$
est non singulier et arbitrairement proche d'un feuilletage par cercles
m\'eridiens.
Le feuilletage caract\'eristique de $T$ pour $\xi_+^t$
est lui constant et dirig\'e par $X$.
Il poss\`ede une certaine pente $p$ (bien d\'efinie modulo $\pi$ par
le choix d'une longitude fixe $l$ de $T$), diff\'erente de celle du m\'eridien
(qui vaudra par d\'efinition $0$)
car sinon on obtiendrait un disque vrill\'e.
Cette pente est la pente
de n'importe quelle courbe int\'egrale de son feuilletage caract\'eristique,
car celui-ci ne poss\`ede que des singularit\'es positives.

On note  $p_+^{\pm 1,t}$  les pentes de $\xi_+^t T\times \{\pm 1\}$.
Si $t$ est assez grand, elles sont non nulles (car sinon $\xi_+$
serait vrill\'ee), d\'ependent contin\^ument de $t$ et tendent vers $0$ (modulo $\pi$)
lorsque $t$ tend vers $+\infty$. On peut alors
\'egalement les supposer irrationnelles.

Si $t$ est assez grand, les pentes $p_1^{1,t}$ et $p_1^{-1,t}$
sont proches de $0$ et en particulier, elles d\'elimitent un secteur
angulaire de petite ouverture qui ne contient pas $p$.
Si on parcourt le cercle dans le sens trigonom\'etrique
pour joindre $p_1^{1,t}$ \`a $p_1^{-1,t}$ en passant par $p$,
on est alors assur\'e de traverser  un intervalle $I(t)$ de
taille minor\'ee par $\pi -\epsilon (t)$, o\`u  $\epsilon (t)$
tend vers $0$ quand $t$ tend vers $+\infty$.

La restriction de $\xi_+^t$ a $T\times [-1,1 ]$ est
universellement tendue~: proche de $\F$, elle est transversale \`a une
fibration en intervalles de $T \times [-1,1]$~; son rappel dans le rev\^etement
universel de ce produit se plonge alors de mani\`ere explicite dans la structure
de contact tendue standard de $\R^3$.
D'apr\`es le th\'eor\`eme principal de \cite{Gi2}
qui classifie les structures de contact tendues sur le tore
\'epais, pour chaque pente $p'$ dans l'intervalle $I(t)$,
il existe un tore $T_{p'}$ dans $T \times [-1,1]$,
isotope \`a $T$ et dont le feuilletage caract\'eristique pour $\xi_+^t$
  est non singulier et  de pente
$p'$. En particulier $I(t)$ ne contient pas $0$ et $\pi$ (modulo $2\pi$),
car sinon $\xi_+^t$, et donc $\xi_+$, serait vrill\'ee.
On note $k$ l'ordre de l'\^ame de la composante de Reeb dans
$H_1 (V,\Z )$. Pour tout $n\in \Z$ premier avec $k$, si $t$
est assez grand, toute courbe de pente
$\frac {k}{n}$ est r\'ealis\'ee  par une
courbe non singuli\`ere $\gamma_{\frac{k}{n}}$
du feuilletage caract\'eristique pour $\xi_+^t$
d'un tore $T_{\frac{k}{n}}$.
Soit $n\in \Z$ premier avec $k$. Pour tout $l\in \Z$,
la courbe $\gamma_{\frac{k}{n+lk}}$ est obtenue comme image
de $\gamma_{\frac{k}{n}}$ par la puissance $l$-i\`eme d'un
twist de Dehn autour du m\'eridien. Comme courbe lisse, elle
est isotope \`a $\gamma_{\frac{k}{l}}$.
Elle  borde donc
une surface compacte et orient\'ee de genre $g$ ind\'ependant de $l$.
Son invariant de Thurston-Bennequin d\'epend affinement de $l$~:
la diff\'erence entre $tb (\gamma_{\frac{k}{n+lk}} )$ et
$tb (\gamma_{\frac{k}{n+(l+1)k}} )$ vaut $k$,
et donc l'in\'egalit\'e de Bennequin
$$tb (\gamma_{\frac{k}{n+lk}} )\leq 2g-1$$ sera viol\'ee
pour un choix de $l$ judicieux.

On en d\'eduit par contradiction que l'holonomie est l'identit\'e au moins d'un c\^ot\'e
de $T$.
Cette argument montre \'egalement que si $U\simeq T^2 \times [0,1]$
est une sous-vari\'et\'e de $V$ feuillet\'ee par des tores
$(T^2 \times \{ s\} )_{s\in [0,1]}$
de $\F$, alors l'holonomie de $\F$ est l'identit\'e pr\`es de
$T^2 \times \{ 0\}$ ou de $T^2 \times \{ 1\}$.  Comme on est parti
d'un tore $T$ bordant une composante de Reeb, et en particulier un tore
solide, on en d\'eduit
que
$V$ est diff\'eomorphe \`a $S^1 \times \R^2$, ce qui donne une contradiction.
\end{proof}

\begin{remarque} Si $(\xi_+ ,\xi_- )$ est une paire positive de
structures de contact tendues, l'in\'egalit\'e de Bennequin montre,
m\^eme lorsque $\lambda$ n'est pas uniquement int\'egrable,
qu'aucune courbe transversale \`a $\lambda$ n'est le bord d'un disque
plong\'e dans $V$.
\end{remarque}
Sous les conditions du th\'eor\`eme~\ref{theoreme : reeb},
le feuilletage $\F$ peut-il poss\'eder des composantes de Reeb~?

L'hypoth\`ese de positivit\'e faite sur  la paire n'est pas tr\`es
agr\'eable \`a manipuler.
On pourrait essayer de la remplacer par une condition d'homotopie~:
une paire de structures de contact $(\xi_+ ,\xi_- )$ tendues
qui sont homotopes comme champs de plans est-elle toujours d\'eformable
en une paire positive~? Autrement dit, si $\xi_+$ et $\xi_-$
sont tendues et homotopes, peut-on toujours \'eliminer leurs contacts
n\'egatifs par une isotopie de $\xi_+$~?
Par ailleurs, il y a des paires positives de structures de contact, m\^eme
transversales, pour lesquelles $\lambda$ n'a pas de feuilletage int\'egral
global~\cite{ET, Mi1,Mi2}.
La notion de paire positive de structures de contact tendues semble
g\'en\'eraliser strictement celle de feuilletage sans
composante de Reeb. Quelles sont les vari\'et\'es qui portent
une telle paire~? Parmi celles-ci, y en a-t-il une qui ne porte
pas de feuilletage sans composante de Reeb~?

Ces commentaires et le  th\'eor\`eme~\ref{theoreme : reeb} nous incitent
\`a formuler la conjecture suivante.

\begin{conjecture} Si $V$ est une vari\'et\'e close qui porte une paire
de structures tendues homotopes, alors le rev\^etement
universel de $V$ est conjugu\'e \`a $\R^3$.
\end{conjecture}

Pour ce qui concerne les structures de contact universellement
tendues sur les vari\'et\'es irr\'eductibles, la preuve de  cette conjecture r\'esulte du fait que la sph\`ere
$S^3$ ainsi que le produit $S^2 \times \R$ ne portent
pas de structures de contact tendues positives et n\'egatives
homotopes~\cite{El}. L'utilisation du th\'eor\`eme
d'uniformisation d\'emontr\'e par  Perelman~\cite{Pe}
permet d'exclure toute autre possibilit\'e que $\R^3$.

Soit $V$ une vari\'et\'e compacte de bord non vide
portant une paire positive $(\xi_+ ,\xi_- )$.
On dit que le bord de $V$ est {\it adapt\'e}
\`a la paire $(\xi_+ ,\xi_- )$ si $\Delta =\Delta_+$ est
transversale \`a $\partial V$ et si $X$ est transversal
\`a $\partial V$ et  sort de $V$
en dehors des points de $\Delta \cap \partial V$. 
C'est par exemple le cas lorsque  $\partial V$ est compos\'e de 
tores transversaux \`a $X$ (et n'intersectant pas $\Delta$) ou 
de sph\`eres intersectant $\Delta$ en deux foyers, comme dans le
th\'eor\`eme~\ref{t:sphere}. Dans la situation o\`u le bord de $V$
est adapt\'e \`a une paire positive $(\xi_+ ,\xi_- )$, 
on obtient, comme dans le cas clos, l'existence d'un champ de
plans $\lambda =\lambda^{+\infty}$, continu et localement int\'egrable.
Les th\'eor\`emes pr\'ec\'edents s'\'etendent ainsi sans difficult\'e 
aux vari\'et\'es compactes dont le bord
est adapt\'e \`a une paire $(\xi_+ ,\xi_- )$.

On donne maintenant  un crit\`ere qui permet d'assurer
que le champ de plans $\lambda$ est uniquement int\'egrable.
Une paire positive de structures de contact est dite
{\it transitive}
si, pour tout $x\in V$, il existe $x_{-\infty} \in \Delta_+$ tel que $X^t (x)$
converge vers $x_{-\infty}$ lorsque $t$ tend vers $-\infty$.
Cette d\'efinition pr\'esente surtout un int\'er\^et dans le cas 
o\`u $V$ poss\`ede un bord non vide adapt\'e \`a $(\xi_+ ,\xi_- )$. 

\begin{theoreme}\label{theoreme : transitif}
Soit $V$ une vari\'et\'e compacte dont le bord est adapt\'e \`a
une paire normale et positive $(\xi_+ ,\xi_- )$.
Si la paire $(\xi_+ ,\xi_- )$ est transitive et si le
champ $X$  poss\`ede la propri\'et\'e
g\'en\'erique de n'avoir ni liaison entre les s\'eparatrices d'une 
naissance-mort et d'une selle,
ni suite de deux connections cons\'ecutives entre selles,  alors
le champ de plans $\lambda$ est uniquement
int\'egrable.
\end{theoreme}
\begin{proof}[D\'emonstration.]  On observe
que $\lambda$ est uniquement int\'egrable en un point
$x$ de $V\setminus \Delta_+$
si et seulement si le champ de droites $\lambda \cap
\nu$ est uniquement int\'egrable, o\`u $\nu$ est
le plan tangent \`a un germe de disque en $x$ transversal \`a
$X$. (Pour obtenir le feuilletage int\'egral de
$\lambda$, on pousse celui de $\lambda \cap \nu$ par le flot de $X$.)
Si un point $x$ de $V$ est dans le bassin d'attraction d'un foyer,
on a vu dans la preuve
de la proposition~\ref{proposition : egalite}
que le champ de plans $\lambda$ est $C^\infty$ au voisinage de
$x$ et donc uniquement int\'egrable.

On suppose \`a pr\'esent que $x$ est sur la s\'eparatrice instable d'une
selle $h$. Par hypoth\`ese de g\'en\'ericit\'e, on est s\^ur qu'au moins une des
deux s\'eparatrices stable de $h$ provient d'un foyer $e$.

On \'etudie d'abord le cas o\`u les deux s\'eparatrices
stables de $h$ proviennent de foyers $e$ et $e'$.
On fixe un petit disque $D$ passant par $x$ et transversal
\`a $X$.
Si $\lambda$ n'est pas uniquement int\'egrable en $x$,
alors le champ de droites orient\'ees  $\lambda \cap TD$, dirig\'e
par un vecteur $Y$,  non plus.
On peut donc trouver un petit arc $\gamma$  inclus dans $D$ et  transversal
\`a $Y$ dont tous les points peuvent \^etre obtenus  en suivant $Y$
\`a partir de $x$
pendant un petit temps positif.
Quitte \`a prendre $\gamma$ assez proche de $x$,
on obtient \'egalement que tous les points de $\gamma$
proviennent de foyers proches de $e$ situ\'es sur une m\^eme
branche de $\Delta_+$ en suivant le flot de $X$.
En particulier, $\lambda$ est $C^\infty$ au voisinage de
$\gamma$ et s'int\`egre uniquement  sur un voisinage $U$
de $\gamma$ en un feuilletage en disques transversal \`a
$\gamma$.

Si  $y$ et $z$ sont deux points assez proches situ\'es sur
$\gamma$, ils doivent  \^etre dans le bassin d'attraction de
deux foyers distincts $e_0$ et $e_1$ proches de $e$ (et situ\'es
sur la m\^eme composante de $\Delta_+$), sinon $y$ et
$z$ seraient sur une m\^eme feuille du feuilletage int\'egral
de  $\lambda \vert_U$.
De ces deux foyers sont issues les s\'eparatrices stables
de deux selles $h_0$ et $h_1$ proches de $h$ (et situ\'ees
sur la m\^eme composante de $\Delta_+$).
Mais alors les feuilles int\'egrales  de $\lambda$
passant par $y$ et $z$, qui sont bien d\'efinies dans
les bassins d'attractions de $e_0$ et $e_1$,
sont asymptotes aux s\'eparatrices instables
de $h_0$ et de $h_1$, ce qui exclut qu'elles puissent
\^etre asymptotes toutes les deux \`a la s\'eparatrice de
$h$.

Si on r\'ep\`ete cet argument pour les petits
temps n\'egatifs de $Y$, on obtient que $Y$
est uniquement int\'egrable au voisinage de $x$,
et donc aussi $\lambda$ (on rappelle que l'unique
feuilletage int\'egral de $\lambda$ est obtenu en poussant
le feuilletage int\'egral de $Y$ par le flot de $X$ qui est
$C^\infty$).

Dans le cas o\`u une s\'eparatrice stable de $h$
provient d'un foyer et l'autre d'une selle $h'$,
on est s\^ur que (par g\'en\'ericit\'e) les deux s\'eparatrices
stables de $h'$ proviennent d'un foyer.
On peut alors appliquer le raisonnement pr\'ec\'edent
le long des deux s\'eparatrices instables de $h'$
pour obtenir que $\lambda$ y est uniquement int\'egrable,
et donc \'egalement au voisinage de celles-ci.

Reste le cas o\`u la deuxi\`eme s\'eparatrice stable
de $h$ provient d'une naissance-mort $n$, g\'en\'eriquement
de l'int\'erieur de la partie foyer. Si $\lambda$ n'\'etait pas
uniquement int\'egrable pr\`es de $x$, on trouverait
une petite transversale $\gamma$ dont tous les points peuvent
\^etre atteints depuis $x$ par des petits arcs tangents \`a $Y$.
Un petit sous-arc $\gamma'$ de $\gamma$ serait constitu\'e de points
issus d'un petit arc de foyers (soient proches de $n$, soient
proches du foyer dont provient g\'en\'eriquement la s\'eparatrice
stable de $n$). De ces foyers sont issues
les s\'eparatrices stables de selles proches de $h$.
On raisonne alors comme
dans le premier cas pour conclure.
Le champ $\lambda$ est donc uniquement int\'egrable
au voisinage des s\'eparatrices stables de $h$.

L'argument d\'evelopp\'e dans le dernier sous-cas fournit
\'egalement que
si l'orbite d'un point $x$ provient
de l'int\'erieur de la partie foyer d'une naissance-mort $n$,
alors $\lambda$ est uniquement int\'egrable au voisinage de
$x$.
Le cas o\`u l'orbite de $x$ est issue du 
bord de la partie foyer d'une naissance-mort
se traite  de mani\`ere identique \`a celle des s\'eparatrices
instables des selles. (On est  dans la situation g\'en\'erique
o\`u la s\'eparatrice stable de la naissance-mort provient
d'un foyer.)
Enfin, le cas o\`u $x \in \Delta_+$ (essentiellement si
$x$ est une selle ou une naissance-mort, puisque celui des foyers
est d\'ej\`a trait\'e) s'\'etudie par des m\'ethodes
similaires.
\end{proof}

Le r\'esultat suivant donne un moyen, par isotopie de $\xi_+$ et
$\xi_-$, de modifier l'ensemble des contacts positifs $\Delta_+$.
Il pourrait permettre dans certains cas de rendre une paire
transitive. C'est une version du lemme de cr\'eation/\'elimination
de singularit\'es de \cite{Gi1} pour les paires.

\begin{proposition}\label{proposition : elimination}
Soient $(\xi_+ ,\xi_- )$ une paire de structures de contact,
$\gamma$ une courbe ferm\'ee positivement transversale \`a $\xi_+$
et $\xi_-$. On suppose que $\gamma \cap \Delta =\emptyset$ et
qu'un voisinage $U\simeq \gamma \times D^2$ de $\gamma \simeq
\gamma \times \{ 0\}$ est tel que $\xi_+ \cap \xi_- =\xi_+ \{ * \}
\times D^2 =\xi_- \{ *\} \times D^2$. Alors, il existe une
isotopie de $\xi_+$ et une isotopie de $\xi_-$ \`a supports dans
$U$ menant \`a une paire $(\xi_+' ,\xi_-' )$ telle que $\xi_+'
\cap \xi_-' =\xi_+'  \{ * \} \times D^2 =\xi_-' \{ *\} \times D^2$
et que les feuilletages caract\'eristiques $\xi_+' \{ * \} \times
D^2 $ et $\xi_-' \{ * \} \times D^2 $ soient tous les deux \'egaux
\`a un m\^eme feuilletage pr\'esentant une selle et un foyer
positifs en position d'\'elimination.
\end{proposition}
\begin{remarque} Dans cette proposition, le lieu des contacts
positifs est modifi\'e par l'ajout de deux composantes \guil{parall\`eles},
l'une form\'ee de la r\'eunion des selles, l'autre de la r\'eunion des
foyers inclus dans $U$.
\end{remarque}

\begin{proof}
Soit $D^2_\epsilon \subset D^2$.
On dessine sur $D^2_\epsilon$ un feuilletage de divergence
positive constitu\'e d'une selle et d'un foyer positif
en position d'\'elimination.
Il est d\'efini comme noyau d'une $1$-forme $\beta$ avec $d\beta >0$.
On munit $S^1 \times D^2_\epsilon$ de la paire de structures
de contact $\zeta_\pm$ d'\'equations $dt\pm \epsilon' \beta =0$, $t\in S^1$.
On d\'epose cette paire sur un voisinage tubulaire $\epsilon$-petit
$U_\epsilon \simeq \gamma \times D_\epsilon$ de $\gamma$,
$D_\epsilon \subset D$,
identifi\'e \`a $S^1 \times D^2_\epsilon$.
On fait en sorte $\zeta_+ \cap \zeta_- =\xi_+ \cap \xi_-$ le long
de $\partial U_\epsilon$. Il est alors possible~\cite{Gi1}, pour $\epsilon'$
et $\epsilon$ assez petits, de raccorder $\zeta_+$ \`a $\xi_+$ et $\zeta_-$ \`a $\xi_-$
par une paire de structures transversales \`a $\frac{\partial}{\partial t}$
et \`a $T(\{ * \} \times D^2 )$
sur $U\setminus U_\epsilon$ --~et donc transversales entre elles~--
et dont l'intersection est dans $T(\{ * \} \times D^2 )$.
La paire obtenue par recollement est la paire $(\xi_+' ,\xi_-' )$ recherch\'ee.
Les structures $\xi_\pm'$ et $\xi_\pm$ sont isotopes comme dans~\cite{Gi1}.
\end{proof}

Au lieu de s'int\'eresser aux paires positives de structures de contact
tendues, on peut s'int\'eresser aux paires pour lesquelles
deux points quelconques de $V$ peuvent \^etre joints
par un arc positivement transversal \`a $\xi_+$ et $\xi_-$
(et donc en particulier par un arc positivement transversal
\`a $\lambda$, ce qui mime la d\'efinition de feuilletage tendu).
Celles-ci sont automatiquement positives.
On parlera alors de paire {\it fortement tendue}.
Dans ce cas, on peut simplifier $\Delta_+$.

\begin{proposition}\label{proposition : transversal} Si $(\xi_+ ,\xi_- )$ est une paire
normale et fortement tendue, alors il existe une isotopie de
$\xi_+$ en une structure $\xi'_+$  telle que la paire $(\xi'_+
,\xi_- )$ soit normale et  fortement tendue et que le lieu des
contacts $\Delta_+$ soit transversal \`a $\xi'_+$ et $\xi_-$ (en
particulier, il n'y a pas de naissance-mort).
\end{proposition}
\begin{proof}[D\'emonstration.] C'est le cas relatif de la preuve
de la
proposition~\ref{proposition : elimination}.
On relie les contacts quadratiques par une collection d'arcs positivements
transversaux \`a $\xi_\pm$ qui \'evitent $\Delta_+$.
Pr\`es des contacts quadratiques, on prend soin que ces arcs soient
situ\'es dans le demi-espace, d\'elimit\'e par $\xi_+ =\xi_-$,
qui ne contient pas $\Delta_+$. En particulier, la direction de $X$
le long de ces arcs converge au bord vers la direction de la demi-s\'eparatrice
de type selle des naissances-morts.
On modifie alors $\xi_+$ et $\xi_-$ le long de chacun de ces arcs pour
faire appara\^ \i tre deux familles de contacts parall\`eles,
l'une constitu\'ee
de foyers, l'autre de selles, qui permettent de faire dispara\^ \i tre
les naissances-morts.
\end{proof}

\begin{remarque} Dans le cas g\'en\'eral, si $\Delta$
est sans contact quadratique avec $\xi_\pm$, on peut d\'efinir
sur $V\setminus \Delta $ deux champs de plans transversaux,
qui se comportent le long de $\Delta$
\`a la mani\`ere du plan tangent \`a un feuilletage d'\'energie
finie pr\`es de son lieu singulier~\cite{HWZ}.
\end{remarque}

\begin{exemple} {\rm Si $(K,\theta )$ est une d\'ecomposition
de $V$ en {\it livre ouvert}, la construction de Thurston et
Winkelnkemper \cite{TW}, syst\'ematis\'ee par Giroux
dans~\cite{Gi3}, fournit une paire $(\xi_+ ,\xi_- )$ de structures
de contact {\it port\'ees} par $(K,\theta )$  pour laquelle
$\Delta_- =K$ (et est constitu\'e de foyers) et dont la
restriction \`a $V\setminus N(K)$ est fortement tendue
(voir~\cite{Gi3} pour une description de ces notions). \`A l'aide
de la proposition~\ref{proposition : transversal}, on peut
\'egalement en obtenir une (normale) pour laquelle  de surcro\^ \i
t $\Delta_+$ est transversale \`a $\xi_+$, ainsi qu'aux pages de
$(K,\theta )$.}
\end{exemple}

Si $\lambda$ est uniquement int\'egrable, et
si $(\xi_+ ,\xi_- )$ est normale et fortement tendue, alors le
feuilletage int\'egral $\F$ de $\lambda$ est tendu.
Comme les structures $\xi_\pm^t$ convergent vers $\F$,
elles sont tendues et aussi $\xi_+$ et $\xi_-$.

R\'eciproquement, si $\xi_+$ et $\xi_-$ sont tendues, forment
une paire normale 
et si $\lambda$ est uniquement int\'egrable, alors $\lambda$ est
sans composante de Reeb homologue \`a z\'ero.
Dans le cas o\`u la vari\'et\'e $V$ est une sph\`ere d'homologie rationnelle
atoro\"\i dale, $\F$ est tendu.
On peut alors montrer que
la paire $(\xi_+ ,\xi_- )$ peut \^etre d\'eform\'ee
(par une isotopie de $\xi_+$)
en une  paire fortement tendue.

\begin{proposition}\label{proposition : fortement}
 Soit $(\xi_+ ,\xi_- )$ une paire normale positive
de structures tendues sur une sph\`ere d'homologie rationnelle
atoro\"\i dale. Si $\lambda$ est uniquement int\'egrable, alors
$(\xi_+ ,\xi_- )$ est d\'eformable, par une isotopie de $\xi_+$,
en une paire normale et fortement tendue.
\end{proposition}
\begin{proof}[D\'emonstration.]
Comme on l'a d\'ej\`a signal\'e plus haut, sous les hypoth\`eses
de la proposition~\ref{proposition : fortement}, le th\'eor\`eme~\ref{theoreme : reeb}
implique que le feuilletage int\'egral
de $\lambda$ est sans composante de Reeb. De plus, comme la vari\'et\'e
ambiante est atoro\"\i dale, il est sans composante de Reeb g\'en\'eralis\'ee,
c'est-\`a-dire qu'il est tendu.

Il suffit de montrer que pour $t$ assez grand,
la paire $(\xi_+^t ,\xi_-^t )$ est fortement tendue. Pour cela,
on va utiliser un syst\`eme d'arcs transversal \`a $\lambda$
et utiliser la proximit\'e entre $\lambda$ et $\xi_\pm^t$.

On se donne deux familles de boules $(B_i (1))_{1\leq i\leq n}$ et
$(B_i (2))_{1\leq i\leq n}$ dans $V$ avec les propri\'et\'es suivantes~:
\begin{itemize}
\item $\Int (B_i (1))$ recouvre $V$~;
\item $B_i (1) \subset \Int B_i (2)$~;
\item il existe un diff\'eomorphisme $\phi : B_i (2)\rightarrow B(0,2)\subset
\R^3$,
avec $\phi (B_i (1))=B(0,1)$ et $\phi_* \lambda$ est $\epsilon$-$C^0$-proche
de $\{ dz=0\}$.
\end{itemize}

Sous ces hypoth\`eses, si $\epsilon$ est petit,
on peut joindre tout point de $B_i (1)$ (resp. le p\^ole sud de
$B_i (2)$) au p\^ole nord de
$B_i (2)$ (resp. tout point de $B_i (1)$) par un arc positivement
transversal \`a $\lambda$,  dont l'angle
avec $\lambda$ est en tout point sup\'erieur \`a une certaine
constante  $c$.

En particulier, si $t$ est assez grand ($>t_0$), on peut joindre
tout point de $B_i (1)$ (resp. le p\^ole sud de
$B_i (2)$) au p\^ole nord de
$B_i (2)$ (resp. tout point de $B_i (1)$) par un arc positivement
transversal \`a $\xi_+$ et $\xi_-$.

On se donne alors pour tout couple $(i,j)$, $1\leq i,j\leq n$, deux  arcs qui relient
les p\^oles nord et sud de $B_i (2)$ aux p\^oles sud et nord
de $B_j (2)$ et qui sont positivement transversaux \`a
$\lambda$. L'existence de ces arcs est assur\'ee par le fait
que le feuilletage int\'egral de $\lambda$ est tendu.
\`A nouveau, si $t$ est assez
grand ($>t_1$), chacun de ces arcs est positivement
transversal \`a  $\xi_+^t$ et $\xi_-^t$.
Pour conclure, on remarque que si $t>max (t_0 ,t_1 )$,
et si $x\in B_i (1)$, $y\in B_j (1)$ sont deux
points quelconques de $V$, on peut joindre
$x$ \`a $y$ par un chemin transversal \`a $\xi_+$ et $\xi_-$
en concat\'enant un chemin entre $x$ et le p\^ole nord
de $B_i (2)$, le chemin entre le p\^ole nord de $B_i (2)$
et le p\^ole sud de $B_j (2)$ et un chemin entre
le p\^ole sud de $B_i (2)$ et $y$.
\end{proof}
\begin{remarque} L'isotopie qui d\'eforme $\xi_+$ ne pr\'eserve pas
$\xi_-$ et on ne peut donc pas dire {\it a priori} que la paire
initiale est fortement tendue.
\end{remarque}

\section{Vari\'et\'es \`a bords sph\'eriques}

Dans le cas des vari\'et\'es bord\'ees par une sph\`ere,
se donner une paire positive tendue \guil{convexe au bord}
est tr\`es fortement contraignant pour la vari\'et\'e.
On notera un  analogue frappant de ce r\'esultat
obtenu {\it via} une m\'ethode de remplissage par des disques holomorphes~\cite{EH}.

\begin{theoreme}\label{t:sphere} Soit $V$ une vari\'et\'e compacte connexe
de dimension trois avec $\partial V =S^2$. Si on suppose que $V$
porte une paire positive $(\xi_+ ,\xi_- )$ de structures tendues
telle que $\Delta_+$ rencontre $\partial V$ transversalement en
deux points $S$ et $N$ et que $X$ soit strictement sortant de $V$
le long de $\partial S\setminus (S\cup N)$, alors $V\simeq B^3$.
\end{theoreme}
\begin{proof}[D\'emonstration.] D'apr\`es la proposition~\ref{proposition : generique},
on peut se ramener au cas o\`u la paire $(\xi_+ ,\xi_- )$ est
normale. La condition au bord nous dit que le champ de plans
$\lambda$ est bien d\'efini sur $V$. Il est transversal \`a
$\partial V$ en dehors de $N$ et $S$. La courbe $\Delta_+$
rencontre $\partial V$ aux points  $N$ et $S$ en un foyer de $X$. En
particulier, $\lambda$ est uniquement int\'egrable, lisse, 
et s'int\`egre en un feuilletage en
disques au voisinage de ces p\^oles.

Soit $A$ l'ensemble des points de $V$ qui appartiennent \`a
un disque int\'egral de $\lambda$, de classe nulle dans
$H_2 (V,\partial V,\Z )$ et le long duquel $\lambda$
est uniquement int\'egrable.
On montre que $A$ est ouvert et ferm\'e dans $V$, et donc \'egal
\`a $V$ par connexit\'e de $V$.

\begin{lemme} $A$ est ferm\'e.
\end{lemme}
\begin{proof}[D\'emonstration.]
Soit $x$ un point de $V$, limite d'une suite $(x_n )_{n\in \N}$ de points
de $A$, situ\'es tous \guil{au-dessus} de $x$, par rapport \`a $\lambda (x)$,
c'est-\`a-dire que, si $D_n$ d\'esigne le disque int\'egral
de $\lambda$ passant par $x_n$, $x$ est situ\'e dans
la composante de $V\setminus D_n$  qui  contient  $S$.
On note $D$ l'ensemble des points d'adh\'erence de
la famille $D_n$.

Si $y\in D$ est la limite d'une suite $(y_{\phi (n)} )$ de points
de $D_{\phi (n)}$, alors dans un petit voisinage cylindrique $U=D^2 \times [-1,1]$
de $y =(0,0,0)$, la composante de $D_{\phi (n)}\cap U$ qui contient
$y_{\phi (n)}$ est le graphe d'une fonction $f_{\phi (n)}$
de classe $C^1$ au-dessus de $D^2$,
et qui est strictement positive en $(0,0)$ car $x$ est au-dessous des
disques $D_n$, et donc aussi $y$.
D'apr\`es le th\'eor\`eme d'Ascoli, une sous-suite de la suite
$(f_{\phi (n)} )$ converge uniform\'ement vers une fonction $f$.
C'est alors aussi le cas de la suite elle-m\^eme car les graphes sont
deux \`a deux disjoints.
La suite des plans tangents aux graphes des $f_{\phi (n)}$ converge elle
aussi, car c'est $\lambda \vert_{D_{\phi (n)}}$, et on en d\'eduit 
que le graphe de $f$ est un disque int\'egral
de $\lambda$, inclus dans $D$.
On remarque maintenant que la composante connexe de $D\cap U$ qui contient
$y$ est \'egale au graphe de $f$. Ceci r\'esulte du fait que
les disques $D_n$ sont deux \`a deux disjoints,
car $\lambda$ est uniquement int\'egrable le long de $D_n$, et tous situ\'es
au-dessus du graphe de $f$.

En conclusion, on obtient que $D$ est une surface int\'egrale de $\lambda$.
Un argument de  Haefliger~\cite{Ha} fournit que  $D$ est compacte,
car limite de feuilles compactes. Comme elle est limite de
disques, on obtient alors \'egalement que $D$ est un disque de classe $C^1$.
Son feuilletage caract\'eristique $\lambda D =\xi_+ D=\xi_- D$
dirig\'e par $X$
est sortant le long du bord. Comme
$\xi_\pm$ sont tendues, il ne poss\`ede pas de cycle limite et
toutes les feuilles sont issues d'une singularit\'e
de $\xi_\pm D=\Delta_+ \cap D$. C'est m\^eme le cas de toute orbite
passant au voisinage de $D$, comme dans le lemme~\ref{lemme : ouvert}.
Le plan $\lambda$ est donc
uniquement int\'egrable le long de $D$ d'apr\`es le
th\'eor\`eme~\ref{theoreme : transitif}.
On raisonne de m\^eme lorsque $x$ est situ\'e \guil{au-dessus} des disques
de la famille $(D_n )$, ce qui permet de conclure dans tous les cas.
\end{proof}

\begin{lemme}\label{lemme : ouvert} $A$ est ouvert.
\end{lemme}
\begin{proof}[D\'emonstration.]

Si $x\in A$ est situ\'e sur un disque int\'egral  $D$ de $\lambda$, toutes
les feuilles de $\xi_\pm  D$ proviennent d'une singularit\'e
de $\xi_\pm D$, c'est-\`a-dire d'un  point de $D\cap \Delta_+$
(comme $\xi_\pm$ sont tendues, $\xi_\pm D$ ne contient pas de cycle limite
et par ailleurs $\xi_\pm D$ est transversal \`a $\partial D$).

Si $\xi_\pm D$ est de type Morse-Smale, toutes les
orbites passant par un point proche de $D$ sont issues de $\Delta_+$
(et m\^eme de foyers ou de selles de $\Delta_+$).
Ceci implique, {\it via} le th\'eor\`eme~\ref{theoreme : transitif},
que $\lambda$ est uniquement
int\'egrable dans un voisinage de $D$.
Un voisinage de $D$ est donc feuillet\'e par des disques
int\'egraux de $\lambda$ d'apr\`es le th\'eor\`eme de stabilit\'e de Reeb.

Si $\xi_\pm D$ poss\`ede une naissance-mort $n$, elle est g\'en\'eriquement
unique sur $D$. Celle-ci est situ\'ee sur une branche $\Delta$
de $\Delta_+$.
La s\'eparatrice stable de $n$ provient d'un foyer $e$ de $\xi_\pm D$
situ\'e sur une composante $\Delta'$ de $\Delta_+$.
D'un c\^ot\'e de $D$, les orbites proches de celles
qui proviennent de $n$ sont toutes issues d'un voisinage de $n$
dans  $\Delta$, et
de l'autre c\^ot\'e de $D$ elles sont toutes issues
d'un voisinage de $e$ dans $\Delta'$.
Comme pr\'ec\'edemment, on conclut en remarquant que $X$ est
transitif au voisinage de $D$.

La m\^eme propri\'et\'e de transitivit\'e au voisinage de $D$
est v\'erifi\'ee si $\xi_\pm D$ pr\'esente une connexion
entre deux selles.
\end{proof}
Ceci montre que $V$ est feuillet\'ee par des disques et
termine la preuve du th\'eor\`eme~\ref{t:sphere}.
 \end{proof}
{\bf Probl\`eme.} Peut-on montrer que si $V$ porte
une paire positive de structures de contact tendues,
alors elle est irr\'eductible~?

$ $\\
Universit\'e de Nantes, Laboratoire de math\'ematiques Jean Leray, UMR 6629 du CNRS, 44322 Nantes, France\\
Vincent.Colin@math.univ-nantes.fr\\
$ $\\
Universidade Federal Fluminense, Niter\'oi, Rio de Janeiro, Brazil\\
firmo@mat.uff.br

\begin{thebibliography}{convex}


\bibitem[Be]{Be}
D.\ Bennequin, {\it Entrelacements et \'equations de Pfaff},
Ast\'erisque {\bf 107-108} (1983), 87--161.


\bibitem[Co1]{Co2} V.\ Colin, {\it Structures de contact tendues sur
les vari\'et\'es toro\"\i dales et approximation de feuilletages
sans composante de Reeb}, Topology {\bf 41} (2002), 1017--1029.

\bibitem[CGH1]{CGH1}
V.\ Colin, E.\ Giroux, and K.\ Honda, {\it On the coarse
classification of tight contact structures,}  Proceedings
of Symposia in Pure Mathematics, {\bf 71}
(2003),  109--120.


\bibitem[El]{El}
Y.\  Eliashberg, \textit{Contact 3-manifolds twenty years since J.\
Martinet's work}, Ann.\ Inst.\ Fourier (Grenoble) \textbf{42} (1992), 165--192.

\bibitem[EH]{EH}
Y.\ Eliashberg et H.\ Hofer, {\it A Hamiltonian characterization of the three-ball},
Journal of Differential and Integral Equations, {\bf 7} No.5 (1994), 1303--1324.


\bibitem[ET]{ET}
Y.\ Eliashberg and W.\ Thurston, \textit{Confoliations}, University Lecture
Series \textbf{13}, Amer.\ Math.\ Soc., Providence (1998).


\bibitem[Gi1]{Gi1}
E.\ Giroux, \textit{Convexit\'e en topologie de contact},
Comment.\ Math.\ Helv.\ \textbf{66} (1991), 637--677.


\bibitem[Gi2]{Gi2}
E.\ Giroux, {\it Structures de contact en dimension trois et
bifurcations des feuilletages de surfaces}, Invent.\ Math.\ {\bf
141} (2000), 615--689.

\bibitem[Gi3]{Gi3}
E.\ Giroux,  {\it G\'eom\'etrie de contact~: de la dimension trois
vers les dimensions sup\'erieures}, Proceedings of the ICM, Vol. II
(Beijing 2002), Higher Ed. Press, Beijing (2002), 405--414.


\bibitem[Gr]{Gr}
J.\ Gray, {\it Some global properties of contact structures}, Ann. of Math.
{\bf 69} (1959), 421--450.

\bibitem[Ha]{Ha}
A.\ Haefliger, {\it Vari\'et\'es feuillet\'ees}, Ann. Scuola Norm. Sup.
Pisa (3) {\bf 16} (1962), 367--397.

\bibitem[Har]{Har} Hartman, {\it Ordinary differential equations},
John Wiley and Sons, Inc. us New York, NY (1964).

\bibitem[HWZ]{HWZ}
H.\ Hofer, K.\ Wyzocki, E.\ Zehnder, {\it Finite energy foliations of tight three spheres and hamiltonian dynamics},
Ann. of Math. {\bf 157} No.1 (2003), 125--255.


\bibitem[Ho1]{Ho1}
K.\ Honda, \textit{On the classification of tight contact structures
I}, Geom.\ Topol.\ {\bf 4} (2000), 309--368.

\bibitem[Ho2]{Ho2}
K.\ Honda, {\it On the classification of tight contact structures
II}, J.\ Differential Geom.\ {\bf 55} (2000), 83--143.


\bibitem[Mi1]{Mi1}
Y.\ Mitsumatsu, \textit{ Anosov flows and non-Stein symplectic manifolds},
Annales de l'institut Fourier, {\bf 45} No.5 (1995), 1407--1421.

\bibitem[Mi2]{Mi2}
Y.\ Mitsumatsu, {\it Foliations and contact structures on
$3$-manifolds}, Foliations: geometry and dynamics (Warsaw, 2000),
World Sci. Publ.,
River Edge, NJ (2002), 75--125.

\bibitem[Pe]{Pe}
G.\ Perelman, {\it The entropy formula for the Ricci flow and
its geometric applications}, arXiv:math/0211159.

\bibitem[Th]{Th}
W.\ Thurston, \textit{ A norm for the homology of
3-manifolds}, Mem.\ Amer.\ Math.\ Soc.\ \textbf{59} (1986), 99--130.

\bibitem[TW]{TW}
W.\ Thurston et H.\ Winkelnkemper, {\it On the existence
of contact forms}, Proc. Amer. Math. Soc. {\bf 52} (1975), 345--347.

\bibitem[Za1]{Za1}
S.\ Zannad, {\it Surfaces branch\'ees en g\'eom\'etrie de contact},
Th\`ese de l'universit\'e de Nantes (2006).

\bibitem[Za2]{Za2}
S.\ Zannad, \textit{A sufficient condition for a branched surface
to fully carry a lamination},
Alg.  Geom. Topol. {\bf 7} (2007), 1599-1632.

\end{thebibliography}
\end{document}